\documentclass[review]{article}

\usepackage{graphicx}

\usepackage{lipsum}
\usepackage{amsfonts}
\usepackage{epstopdf}
\usepackage{algorithmic}

\usepackage{authblk}
\usepackage{bibentry}
\usepackage{caption}
\usepackage{colortbl}
\usepackage{amsmath}
\usepackage{amsfonts}
\usepackage{url}
\usepackage{mathtools}
\usepackage{array}
\usepackage{multirow}
\usepackage{booktabs}
\usepackage{fullpage}
\usepackage{xcolor}

\usepackage{subfig}

\newtheorem{mydef}{Definition}[section]

\title{Computational modelling of cardiac ischaemia using a variable-order fractional Laplacian}

\author[1,2]{Megan E. Farquhar\thanks{me.farquhar@qut.edu.au}}
\author[1,2]{Timothy J. Moroney\thanks{t.moroney@qut.edu.au}}
\author[1,2]{Qianqian Yang\thanks{q.yang@qut.edu.au}}
\author[1,2]{Ian W. Turner\thanks{i.turner@qut.edu.au}}
\author[1,2,3]{Kevin Burrage\thanks{kevin.burrage@qut.edu.au}}

\affil[1]{\footnotesize{School of Mathematical Sciences, Queensland University of Technology, Brisbane, QLD 4000, Australia.}}
\affil[2]{Australian Research Council Centre of Excellence for Mathematical and Statistical Frontiers, Queensland University of Technology, Brisbane, QLD 4000, Australia.}
\affil[3]{Visiting Professor, Department of Computer Science, University of Oxford, Oxford, UK.}

\usepackage{amsopn}
\DeclareMathOperator{\diag}{diag}

\begin{document}

\maketitle

\begin{abstract}
  Heart failure is one of the most common causes of death in the western world. Many heart problems are linked to disturbances in cardiac electrical activity, such as wave re-entry caused by ischaemia. In terms of mathematical modelling, the monodomain equation is widely used to model electrical activity in the heart.  Recently, Bueno-Orovio et al. [J. R. Soc. Interface 11: 20140352, 2014] pioneered the use of a fractional Laplacian operator in the monodomain equation to account for the complex heterogeneous structures in heart tissue. In this work we consider how to extend this approach to apply to hearts with regions of damaged tissue.  This requires the use of a fractional Laplacian operator whose fractional order varies spatially.  We develop efficient numerical methods capable of solving this challenging problem on domains ranging from simple one-dimensional intervals with uniform meshes, through to full three-dimensional geometries on unstructured meshes.  Results are presented for several test problems in one dimension, demonstrating the effects of different fractional orders in regions of healthy and damaged tissue.  Then we showcase some new results for a three-dimensional fractional monodomain equation with a Beeler-Reuter ionic current model on a rabbit heart mesh.  These simulation results are found to exhibit wave re-entry behaviour, brought about only by varying the value of the fractional order in a region representing damaged tissue.
\end{abstract}

\textbf{Keywords} Composite medium, Variable-order fractional Laplacian, Matrix transfer technique, Matrix functions, Monodomain Beeler-Reuter model, Ischaemia, Wave re-entry

\newcommand{\shade}{\cellcolor[gray]{0.8}}
\newcommand{\mdef}{\cellcolor{gray!10}}

\section{Introduction}
Heart failure is the most common cause of death in the western world. An improved understanding of how the heart works may lead to important tools for diagnosis and treatment options for heart problems, or diseased hearts.  To this end, mathematical modelling has an important role to play in improving our understanding of the propagation of electrical signals through the heart.

A standard equation used in cardiac modelling is the monodomain equation \cite{bueno13fractional,sundnes06computing,sundnes06computational,whiteley06an}, a partial differential equation (PDE) given by
\begin{equation}
\frac{\partial v}{\partial t} = -D\left(-\nabla^2\right)v - \frac{1}{C_m}\left(I_{ion} - I_{stim}\right)\,.
\label{eq:mono}
\end{equation}
This equation models the evolution of the potential difference $v$ (also known as the action potential or transmembrane potential) across the membrane that separates the intracellular and extracellular domains.
The parameter $D$ represents an effective
conductivity tensor, the parameter $C_m$ is the membrane capacitance per unit area, $I_{ion}$ is the ionic current across the membrane and $I_{stim}$ is the stimulus current applied to the tissue.

To describe the cellular electrophysiological dynamics the monodomain equation is coupled to a system of ordinary differential equations (ODEs) to model the ionic current.
In this work we consider the Beeler-Reuter ionic current model \cite{beeler77reconstruction}, where the ionic current is given by the sum of the currents driven by Sodium ($I_{Na}$), Potassium ($I_K$ and $I_x$) and Calcium ($I_s$),
\begin{equation}
I_{ion} = I_{Na} + I_{K} + I_x + I_s.
\label{eq:BRode}
\end{equation}
The full equations for each of these currents can be found in Section~\ref{sec:BRmodel}. These currents are controlled by seven gating variables $G=\{m,h,j,d,f,x,c\}$ that are themselves governed by the following ODE system,
\begin{equation}
\begin{aligned}
\frac{dG}{dt} &= \alpha_G(v)(1-G)-\beta_G(v)G, \ \ G = \{m,h,j,d,f,x\}\\
\frac{dc}{dt} &= 0.07(1-c) - I_s\,,
\end{aligned}
\label{eq:BRode2}
\end{equation}
where $I_{stim}$ is an external stimulation current that is applied using some stimulation protocol.

Recent work by Bueno-Orovio et al. \cite{bueno13fractional} suggests that the monodomain model can not fully capture the modulation of impulse propagation caused by structural heterogeneities in heart tissue.  They proposed the \emph{fractional} monodomain model,
\begin{equation}
\frac{\partial v}{\partial t} = -D\left(-\nabla^2\right)^{\alpha/2}v - \frac{1}{C_m}\left(I_{ion} - I_{stim}\right), \ \ \ 1 < \alpha \leq 2
\label{eq:fracMono}
\end{equation}
where the diffusion term now incorporates the fractional Laplacian operator $\left( - \nabla^2\right)^{\alpha/2}$. The fractional power of the Laplacian operator is defined, in the standard way for a self-adjoint operator, through its spectral decomposition.
\begin{mydef}
	\cite{ilic06numerical} Suppose the Laplacian $\left( -\nabla^2\right)$ has a complete set of orthonormal eigenfunctions $\phi_{n}$ corresponding to eigenvalues $\lambda_{n}^2$ on a bounded region $\Omega$
	i.e., $\left(-\nabla^2\right)\phi_{n} = \lambda_{n}^2\phi_{n}$ on $\Omega$; $\mathcal{B}(\phi) = 0$ on $\partial \Omega$, where $\mathcal{B}$ is one of the standard three homogeneous boundary conditions (Dirichlet, Neumann or mixed). Let
	\begin{equation*}
	\mathcal{F}_\alpha = \left\{ f = \sum_{n=1}^\infty  c_{n}\phi_{n}, \ \ c_{n} = \langle f, \phi_{n} \rangle \ \ \left| \ \ \sum_{n=1}^\infty |c_{n}|^2|\lambda_{n}|^{\alpha} < \infty, \ \ 1 < \alpha \leq 2\right. \right\}.
	\end{equation*}
	Then for any $f \in \mathcal{F}_\alpha$, $(-\nabla^2)^{\alpha/2}f$ is defined by
	\begin{equation*}
	(-\nabla^2)^{\alpha/2}f = \sum_{n=1}^\infty c_{n} \lambda_{n}^{\alpha} \phi_{n}.
	\end{equation*}
	\label{def:fracLap}
\end{mydef}

Bueno-Orovio et al. found that a fractional value of $\alpha$, commonly referred to as the fractional \emph{order}, better captured the action potential foot width exhibited in experimental tissue depolarisation results than the standard diffusion model with $\alpha = 2$ (their results suggested that $\alpha = 1.75$ was the appropriate value).  Furthermore, they showed that the concave phase plane trajectories observed during depolarisation could \emph{only} be replicated with a fractional order; standard diffusion produces a linear profile.

\begin{figure}
	\centering
	\subfloat[Diffuse myocardial fibrosis (Ischaemic fibrosis of the myocardium) (Simionescu trichromic staining, ob. x4) : myocardial cells (red) intermingled with collagen-rich fibrosis (blue) that completely replaced the necrotic myocardial cells. Capillaries (with yellow-orange red blood cells) within fibrosis remained from repair by the connective tissue process.]{\includegraphics[width = 0.45\textwidth]{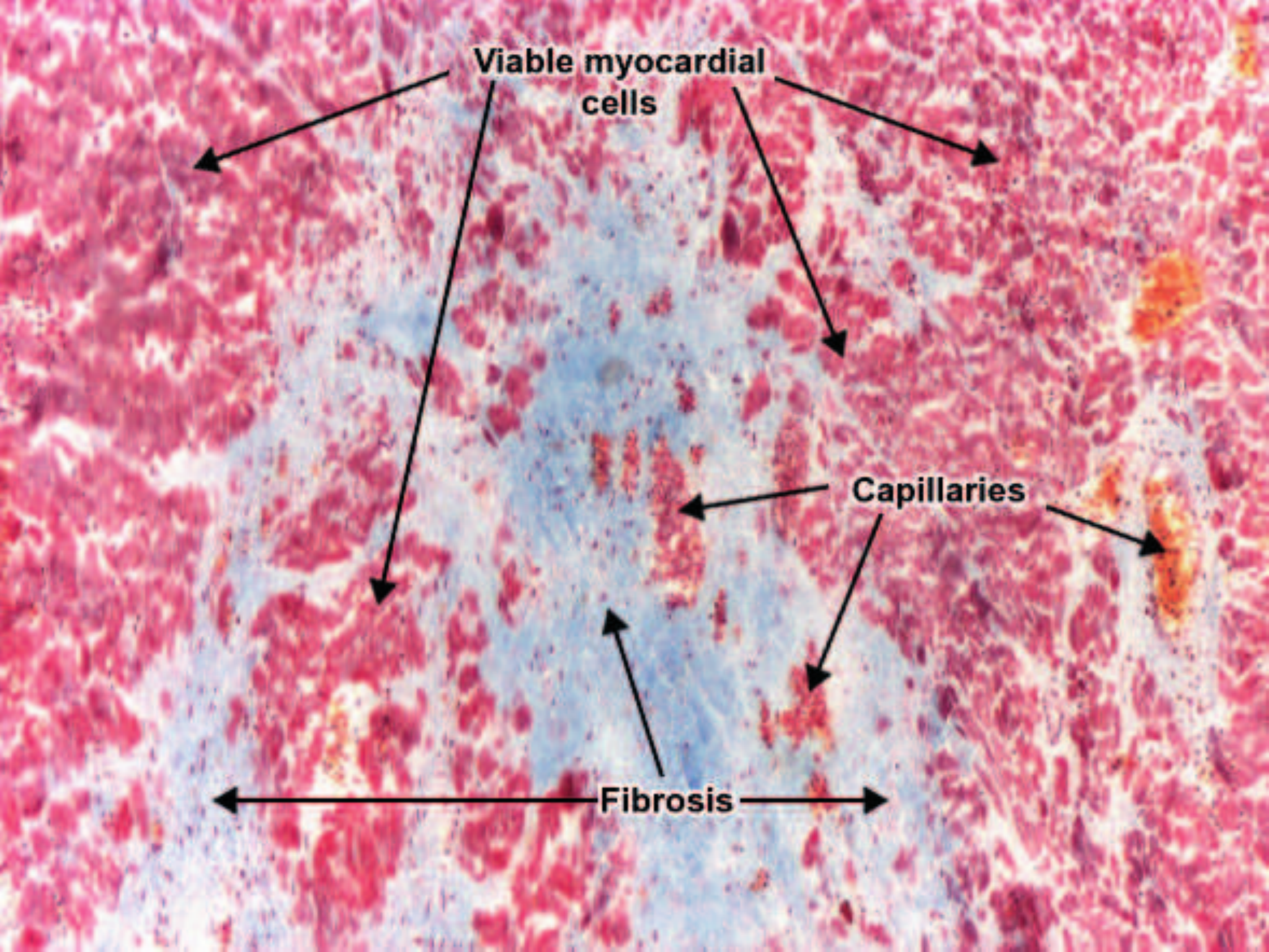}}
	\qquad
	\subfloat[Diffuse myocardial fibrosis (Ischaemic fibrosis of the myocardium) (Simionescu trichromic staining, ob. x10) : viable myocardial cells (red) with nuclei (brown) surrounded by collagen-rich fibrosis (blue). Fibrosis completely replaced the necrotic myocardial cells. Capillaries (with yellow-orange red blood cells) within fibrosis remained from repair by the connective tissue process.]{\includegraphics[width = 0.45\textwidth]{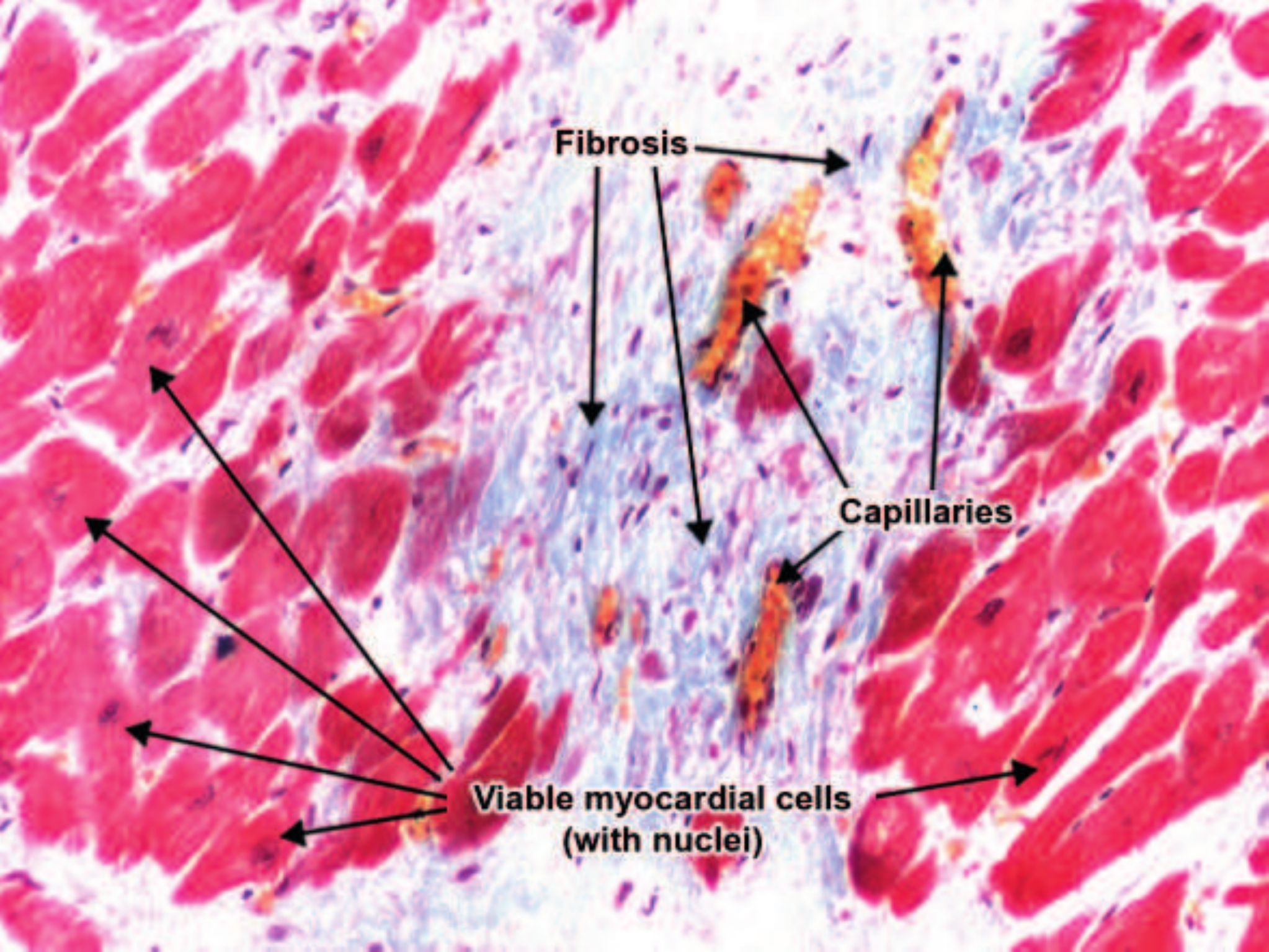}}
	\caption{Images of myocardial fibrosis at the cellular scale.~\cite{atlas10figure}}
	\label{fig:ischaemia4}
\end{figure}

In this work, we are also concerned with the variability in cellular structure that can be present in different regions as a result of damage to the heart.  During a heart attack, regions of cardiac tissue are starved of oxygen, which results in the permanent death of this tissue.  The tissue properties in the damaged (ischaemic) region can be significantly different than the properties in the healthy (normal) tissue region, as seen in Figure~\ref{fig:ischaemia4}. This figure shows cellular level images of heart tissue where ischaemia is present (represented by blue) surrounded by healthy tissue (represented by red), illustrating the difference in cellular structure in the healthy tissue compared to the ischaemic region.

An interesting and important phenomenon that is observed in hearts with ischaemia is wave re-entrant behaviour \cite{dutta16early}. The variation in the excitability and repolarisation properties between the normal and ischaemic regions leads to the establishment of this re-entry behaviour, where the action potential re-enters the system around the ischaemic region, propagating into recovered tissue.

Previous mathematical models have incorporated differences between the normal and ischaemic regions in terms of both the equations and their parameters in the ischaemic region. Dutta et al. \cite{dutta16early} modelled ischaemic re-entry, basing their ventricular membrane kinetics on the ten Tusscher and Panfilov ionic current model \cite{ten06alternans}. They incorporated ischaemic tissue by making changes to the ionic current model in both the ischaemic region and a border region between the ischaemic tissue and the normal tissue. The changes made to the model included increasing the extracellular potassium concentration, decreasing the peak conductance of the fast sodium current, decreasing the L-type calcium current, increasing the resting potential and adding the ATP-sensitive potassium current of Michailova et al. \cite{michailova05modeling}.

Clayton et al. \cite{clayton02reentry} considered the Luo-Rudy ionic current model and incorporated ischaemia by making changes to the model. They individually increased potassium outflow, increased extracellular potassium concentration and decreased the inward sodium current before combining them in a simulation to model moderate and severe ischaemia.

In the context of the fractional model \eqref{eq:fracMono} a new possibility is presented: the fractional order $\alpha$ can be \emph{spatially-varying}, i.e. a function $\alpha(x)$ may be used for the fractional order, \emph{viz.}
\begin{equation}
\frac{\partial v}{\partial t} = -D\left(-\nabla^2\right)^{\alpha(x)/2}v - \frac{1}{C_m}\left(I_{ion} - I_{stim}\right), \ \ \ 1<\alpha(x) \leq 2\,.
\label{eq:varfracmono}
\end{equation}
In fact, in this paper we demonstrate that re-entry behaviour is observed in the simulations of \eqref{eq:varfracmono} by varying \emph{only} the fractional index $\alpha(x)$ between the healthy and ischaemic regions, via a piecewise form
\begin{equation}
\alpha(x) = \left\{ \begin{aligned}
\alpha_1\,, & \ \ \  x \in \Omega_1 \\
\alpha_2\,, & \ \ \  x \in \Omega_2
\end{aligned}\right.
\label{eq:alpha}
\end{equation}
with the rest of the coupled monodomain-ionic current model in \eqref{eq:BRode} and \eqref{eq:BRode2} unchanged.

To perform these simulations, we will introduce a numerical method for solving problems of the general type
\begin{equation}
\frac{\partial u}{\partial t} = - D \left( - \nabla^2 \right)^{\alpha(x)/2}u + g(u)
\label{eq:varpde}
\end{equation}
subject to appropriate boundary conditions.  Before proceeding, we must define precisely what is meant by the variable-order fractional Laplacian $\left(-\nabla^2\right)^{\alpha(x)/2}$.  On \emph{unbounded} domains, this operator can be identified with the variable-order Riesz fractional derivative~\cite{samko13fractional}, given by
\begin{equation}
\mathbb{D}^{\alpha(\cdot)} f(x) = c_{n,\alpha(x)/2}\int_{\mathbb{R}^n} \frac{f(x) - f(x-y)}{|y|^{n + \alpha(x)}} \textrm{d} y
\label{eq:rop}
\end{equation}
with
\begin{equation}
c_{n, \alpha(x)/2} = \frac{2^{\alpha(x)} \Gamma(n/2 + \alpha(x)/2)}{\pi^{n/2} |\Gamma(-\alpha(x)/2)|}.
\end{equation}
For one-dimensional \emph{bounded} domains with homogeneous Dirichlet conditions, it is straightforward to represent the Riesz derivative in terms of two-sided Riemann-Liouville derivatives, where the solution outside the domain is simply set to zero, which is a standard approach in the literature \cite{samko13fractional}.

Dealing with homogeneous \emph{Neumann} boundary conditions, even in one dimension, is not so straightforward. Cusimano \cite{cusimano15fractional} showed how imposing \emph{reflecting} conditions on the boundary leads to a formulation of the variable-order Riesz operator (or rather, its equivalent in terms of two-sided Riemann-Liouville fractional derivatives) that correctly models the Neumann condition in one spatial dimension.  The variable order case where $\alpha(x)$ has the piecewise form \eqref{eq:alpha} was also considered within this framework. Numerically, the variability of the fractional order is naturally incorporated in this approach by forming row $i$ of the (dense) coefficient matrix with the appropriate value of $\alpha$ for the region associated with node $i$.

Our approach is to spatially discretise the variable-order fractional Laplacian operator itself.  Although this poses some new challenges, we will show how the resulting formulation can be made efficient, and by its nature it is applicable generally for Neumann boundary conditions in any number of dimensions.

The standard approach for discretising $\left(-\nabla^2\right)^{\alpha/2}$ (fixed order) is the matrix transfer technique~\cite{farquhar16GPU,burrage12an, yang10novel, yang10numerical}: derive a matrix representation $\mathbf{A}$ of $\left(-\nabla^2\right)$ through finite differences, finite volumes, or similar, then compute with $\mathbf{A}^{\alpha/2}$.

\begin{mydef}
Let $\mathbf{A}$ be a symmetric positive semi-definite matrix representation of the standard Laplacian operator subject to the relevant homogeneous boundary conditions (Dirichlet, Neumann or mixed).  Denote its orthogonal diagonalisation
\begin{equation}
\mathbf{A} = \mathbf{V} \boldsymbol{\Lambda} \mathbf{V}^T
\end{equation}
where
\begin{equation}
\boldsymbol{\Lambda} = \diag(\lambda_1, \lambda_2, \dots, \lambda_N), \ \ \ \ \mathbf{V} = [v_1, v_2, \dots, v_N]
\end{equation}
and $\lambda_i$ is the $i^{th}$ eigenvalue and $v_i$ is the corresponding eigenvector. Under the matrix transfer technique, the operator
$\left(-\nabla^2\right)^{\alpha/2}$ applied to a function $u(x)$ at the $i^{th}$ spatial node $x_i$ is represented by the spatially discretised equation
\begin{equation}
\begin{aligned}
\left.\left(-\nabla^2\right)^{\alpha/2}u(x)\right|_{x = x_i} &\approx \mathbf{e}_i^T (\mathbf{A}^{\alpha/2}\mathbf{u}) \\
&= \sum_j c_j \lambda_j^{\alpha/2} v_{ij},
\end{aligned}
\label{eq:discfraclap}
\end{equation}
where $v_{ij}$ is the element in $\mathbf{V}$ in the $i^{th}$ row, $j^{th}$ column (that is, the $i^{th}$ element of $v_j$)
and
\begin{equation*}
c_j = v_j^T \mathbf{u}.
\end{equation*}
\label{def:discrFracLap}
\end{mydef}

In \cite{farquhar16GPU} we derived an efficient, dense-matrix-free algorithm for calculating the matrix function vector products required in \eqref{eq:discfraclap}.  Dealing only with the sparse matrix representation $\mathbf{A}$, the approach scales well to problems with even millions of spatial nodes.

The generalisation of Definition~\ref{def:discrFracLap} to variable order follows by analogy with \eqref{eq:rop}: the fractional order becomes spatially-dependent with value $\alpha_i \equiv \alpha(x_i)$ at the $i$th spatial node:
\begin{equation}
\left.\left(-\nabla^2\right)^{\alpha(x)/2}u(x)\right|_{x = x_i}= \sum_j c_j \lambda_j^{\alpha_i/2} v_{ij}.
\label{eq:discvarfraclap}
\end{equation}
However, we immediately note an apparent impediment to applying the matrix transfer technique for this variable order case: the right hand side of \eqref{eq:discvarfraclap} is no longer a matrix function, as it was in \eqref{eq:discfraclap}.  Instead, the expression depends explicitly on the diagonalisation of $\mathbf{A}$, via its eigenvalues being raised to a position-dependent power $\alpha_i$.  Nevertheless, despite initial appearances, in the next section we show how the implementation developed in \cite{farquhar16GPU} can be adapted to this new situation.

Furthermore, our formulation is fully extensible to higher dimensions, for complex geometries discretised with unstructured meshes (triangular, tetrahedral, etc.).  We make use of this extensibility in this paper, presenting results for three-dimensional simulations on unstructured tetrahedral meshes in Section~\ref{sec:3dBRres}.  Finally, we note that in the particular case of a one-dimensional domain, with uniform spatial discretisation and Neumann boundary conditions, we recover an equivalent scheme to that of Cusimano et al. \cite{cusimano15fractional}, and comparisons to their work including comments on efficiency will also be made in Section~\ref{sec:nicoleres}.

The remainder of this paper is laid out as follows.  In Sections \ref{sec:vofl} and \ref{sec:numeval} we present our techniques for spatially and temporally discretising the variable-order fractional reaction diffusion equation \eqref{eq:varpde}, and for the efficient numerical evaluation of the required matrix function vector products.  In Section ~\ref{sec:results4} we present results of a number of numerical experiments with variable-order fractional reaction-diffusion equations in one and three dimensions.  These include the aforementioned comparisons with Cusimano et al.'s work and further results simulating the progression of electrical impulses through a `cable' of tissue in one dimension.  Finally, we present simulations through a heart geometry in three dimensions that exhibit re-entry behaviour.

\section{Computing the Variable-order Fractional Laplacian}
\label{sec:vofl}

The discretisation \eqref{eq:discvarfraclap} of the variable-order fractional Laplacian is expressed in terms of the diagonalisation of the matrix representation $\mathbf{A}$.  Using this discretisation directly for computation is unattractive, since diagonalisation is costly in terms of both runtime and of memory storage requirements.  Our idea is to exploit the fact that in the present biological application, the function $\alpha(x)$ has the piecewise form \eqref{eq:alpha}, corresponding to distinct regions of healthy tissue and damaged tissue.

Manipulating Equation~\eqref{eq:discvarfraclap} by adding and subtracting terms involving $\alpha_1$ only, we can obtain the expression for the $i$-th node in terms of matrix functions:
\begin{equation*}
\begin{aligned}
\left.\left(-\nabla^2\right)^{\alpha(x)/2}u(x)\right|_{x = x_i} &\approx \sum_{j=1}^N c_j \lambda_j^{\alpha_i/2} v_{ij}\\
&= \sum_{j=1}^N c_j \lambda_j^{\alpha_i/2} v_{ij} + \sum_{j=1}^N c_j \lambda_j^{\alpha_1/2} v_{ij} - \sum_{j=1}^N c_j \lambda_j^{\alpha_1/2} v_{ij}\\
&= \sum_{j=1}^N c_j \lambda_j^{\alpha_1/2} v_{ij} + \left\{ \begin{aligned}
0& \ \ \ x_i \in \Omega_1 \\
\sum_{j=1}^N c_j (\lambda_j^{\alpha_2/2} - \lambda_j^{\alpha_1/2}) v_{ij} &  \ \ \ x_i \in \Omega_2  \end{aligned}
\right.\\
&= \mathbf{e}_i^T  (\mathbf{A}^{\alpha_1/2}\mathbf{u})
+ \left\{ \begin{aligned}
0& \ \ \ x_i \in \Omega_1 \\
\mathbf{e}_i^T  \left((\mathbf{A}^{\alpha_2/2} - \mathbf{A}^{\alpha_1/2})\mathbf{u}\right) &  \ \ \ x_i \in \Omega_2  \end{aligned}
\right. .
\end{aligned}
\end{equation*}
Supposing that the nodes are ordered first by region, repeating the above for each node in the domain we express the full matrix representation of the variable-order fractional Laplacian operator in terms of matrix functions,
\begin{equation}
\left(-\nabla^2\right)^{\alpha(x)/2}u \Rightarrow \mathbf{A}^{\alpha_1/2}\mathbf{u} +
\mathbf{E}_{\Omega_2}(\mathbf{A}^{\alpha_2/2} - \mathbf{A}^{\alpha_1/2})\mathbf{u}
\label{eq:vopsplit}
\end{equation}
where $\mathbf{E}_{\Omega_i}$ denotes a block matrix with its only nonzero block being the identity in the rows and columns that correspond to nodes $x \in \Omega_i$:
\begin{equation}
\mathbf{E}_{\Omega_1} = \begin{bmatrix}
\mathbf{I}_{N_1} & \mathbf{0}\\
\mathbf{0} & \mathbf{0}
\end{bmatrix},
\mathbf{E}_{\Omega_2} = \begin{bmatrix}
	\mathbf{0} & \mathbf{0}\\
	\mathbf{0} & \mathbf{I}_{N_2}
\end{bmatrix} \textrm{ and } \mathbf{E}_{\Omega_1} + \mathbf{E}_{\Omega_2} = \mathbf{I}_N,
\end{equation}
where $N_i$ is the number of nodes in $\Omega_i$ and $N_1+N_2 = N$.
We note that representation \eqref{eq:vopsplit} could be alternatively written as
\begin{equation}
\left(-\nabla^2\right)^{\alpha(x)/2}u \Rightarrow \mathbf{E}_{\Omega_1}\mathbf{A}^{\alpha_1/2}\mathbf{u} +
\mathbf{E}_{\Omega_2}\mathbf{A}^{\alpha_2/2}\mathbf{u}\,,
\label{eq:badformation}
\end{equation}
however this representation seems less attractive from a numerical point of view because any kind of implicit time differencing cannot be solely written in terms of matrix-function vector products.  Representation \eqref{eq:vopsplit} has the advantage of treating the ``principal'' value $\alpha_1$ (as representative of healthy tissue) in full matrix-function form, while the value $\alpha_2$ corresponding to damaged tissue contributes a correction term.

Using \eqref{eq:vopsplit} to spatially discretise Equation~\eqref{eq:varpde} leads to the following semi-discrete ordinary differential equation system
\begin{equation}
\frac{\textrm{d} \mathbf{u}}{\textrm{d}t} =  -D\mathbf{A}^{\alpha_{1}/2}\mathbf{u} +
\mathbf{E}_{\Omega_{2}}(D\mathbf{A}^{\alpha_{1}/2} - D\mathbf{A}^{\alpha_{2}/2})\mathbf{u}+ g(\mathbf{u})\,.
\end{equation}
A fully implicit backward Euler discretisation with stepsize $\delta t$ is used to complete the discretisation in time, to obtain the fully discrete scheme
\begin{equation}
\mathbf{u}_{n+1} =  f_b(\mathbf{A}) \left( \mathbf{u}_n + \delta t \mathbf{E}_{\Omega_{2}}f_a(\mathbf{A})\mathbf{u}_{n+1}+ \delta t g(\mathbf{u}_{n+1}) \right)\,
\label{eq:BEfp}
\end{equation}
where
\begin{equation*}
f_a(\mathbf{A}) = (D\mathbf{A}^{\alpha_{1}/2} - D\mathbf{A}^{\alpha_{2}/2})\ \ \ \textrm{and} \ \ \ f_b(\mathbf{A}) = \left(\mathbf{I} + D\delta t \mathbf{A}^{\alpha_{1}/2} \right)^{-1} .
\end{equation*}
The fully implicit nature of the scheme is apparent from the appearance of $\mathbf{u}_{n+1}$ on the right hand side of Equation~\eqref{eq:BEfp}. We implement a fixed-point (Picard) iteration to converge the solution with each timestep.  

\section{Numerical Evaluation}
\label{sec:numeval}
The success of the fully implicit scheme described by Equation~\eqref{eq:BEfp} hinges on having an efficient method to evaluate the matrix function vector products in the equation at each time step. In the most general case we will have $\alpha_1 \neq 2$ and $\alpha_2 \neq 2$, meaning that two genuine matrix functions appear in \eqref{eq:BEfp},
\begin{equation*}
f_a(\mathbf{A}) = (D\mathbf{A}^{\alpha_{1}/2} - D\mathbf{A}^{\alpha_{2}/2}) \ \ \ \textrm{and} \ \ \ f_b(\mathbf{A}) = \left(\mathbf{I} + D\delta t \mathbf{A}^{\alpha_{1}/2} \right)^{-1}.
\end{equation*}
We evaluate the matrix function vector products in Equation~\eqref{eq:BEfp} using an approach that combines the contour integral method of Hale et al. for representing matrix functions in terms of families of shifted linear systems \cite{higham08functions} and Krylov subspace methods~\cite{saad03iterative} for solving the linear systems.  We present the key elements of our approach below -- the full details can be found in Farquhar et al. \cite{farquhar16GPU}.

Hale et al.'s ~\cite{trefethen08computing} approach to compute $f(\mathbf{A})\mathbf{b}$ is based on numerically evaluating the contour integral that defines the matrix function vector product (this being one of several, equivalent definitions of a matrix function; see Higham \cite{higham08functions} for alternatives)
\begin{equation}
f(\mathbf{A})\mathbf{b} = \frac{1}{2\pi i} \oint_C f(z) (z\mathbf{I} - \mathbf{A})^{-1}\mathbf{b} \textrm{ d}z,
\end{equation}
where $C$ is a closed contour that lies in the region of analyticity of $f$ and wraps once around the spectrum $\sigma(\mathbf{A})$ in the anticlockwise direction.
Via a clever choice of conformal map, this integral may be recast in a form that yields geometric convergence when evaluated using the midpoint rule \cite{trefethen08computing}.  The transformed integral takes the form
\begin{equation}
f(\mathbf{A})\mathbf{b} = \text{Im} \int_{-K + i K^\prime/2}^{K + i K^\prime/2}  w(\tau)(z(\tau)\mathbf{I} - \mathbf{A})^{-1}\mathbf{b} \ \textrm{d}\tau\,,
\label{eq:int}
\end{equation}
where
\begin{equation}
\label{eq:kdef}
w(\tau)= \frac{-2\sqrt{\lambda_1 \lambda_N}f(z(\tau))\text{cn}(\tau)\text{dn}(\tau)}{\pi k (k^{-1} - \text{sn}(\tau))^2},\qquad k = \frac{\sqrt{\kappa} - 1}{\sqrt{\kappa} + 1}, \qquad \kappa = \frac{\lambda_N}{\lambda_1}\,,
\end{equation}
$\textrm{sn}(\tau) = \textrm{sn}(\tau|k)$, $\textrm{cn}(\tau) = \textrm{cn}(\tau|k)$ and $\textrm{dn}(\tau) = \textrm{dn}(\tau|k)$ are the Jacobi elliptic functions, $\lambda_1$ and $\lambda_N$ are the smallest and largest eigenvalues of $\mathbf{A}$, $K$ and $K^\prime$ are defined by full elliptic integrals with the parameters $k$ and $k^\prime = \sqrt{1-k^2}$, 
\begin{equation*}
K(k) = \int_0^{\pi/2} \frac{\textrm{d}\theta}{\sqrt{1 - k^2 \sin^2(\theta)}},\qquad K^\prime(k) = K(k^\prime)
\end{equation*}
and the contour is given by
\begin{equation}
z(\tau) = \sqrt{\lambda_1\lambda_N}\left(\frac{k^{-1} + \text{sn}(\tau)}{k^{-1} - \text{sn}(\tau)}\right).
\end{equation}
The midpoint rule approximation then takes the form
\begin{equation}
f(\mathbf{A})\mathbf{b} \approx \mathcal{I}_P = \text{Im} \sum_{j=1}^P w(\tau_j)  (z(\tau_j)\mathbf{I} - \mathbf{A})^{-1} \mathbf{b} \,,
\label{eq:contoursol}
\end{equation}
with $P$ the number of quadrature points and $\tau_j = -K + {K(2j-1)}/{P} + ( {K^\prime}/{2}) \ \textrm{i} $ for $j = 1, 2, \dots, P$ and $\textrm{i} = \sqrt{-1}$.

The family of linear system solves required in \eqref{eq:contoursol} involve the shifted matrices $z(\tau_j)\mathbf{I} - \mathbf{A}$, which retain the symmetry and sparsity of $\mathbf{A}$ ($\mathbf{A}$ being, as always, the sparse matrix representation of the Laplacian).  Krylov subspace methods are an attractive option for solving such a family of linear systems, since the Lanczos decomposition is invariant to such shifts.  Hence, only one decomposition, based on the unshifted matrix $\mathbf{A}$, need be computed.

To accelerate the process and reduce the memory required for storing subspace vectors, we apply dual levels of preconditioning to the Lanczos algorithm \cite{farquhar16GPU}.  At the outer level, a deflation preconditioner is applied, which serves to reduce both the number of terms $P$ required in \eqref{eq:contoursol}, and the number of Lanczos iterations required for each term.  This deflation preconditioner uses pre-computed spectral information of $\mathbf{A}$ to deflate its smallest $\ell$ eigenvalues.  Then, using the known relationship between the matrix functions of $\mathbf{A}$ and of the deflated matrix $\hat{\mathbf{A}}$ \cite{farquhar16GPU,yang11novel}, we compute
\begin{equation}
f(\mathbf{A})\mathbf{b} = \mathbf{Q}_\ell f(\mathbf{\Lambda}_\ell) \mathbf{Q}_\ell^T\mathbf{b} + f(\hat{\mathbf{A}})\mathbf{\hat{b}},
\label{eq:preconfunc}
\end{equation}
where $\boldsymbol\Lambda_\ell$ is the diagonal matrix of the smallest $\ell$ eigenvalues, $\mathbf{Q}_\ell$ is the matrix of associated orthonormal eigenvectors and $\mathbf{\hat{b}} = (\mathbf{I} - \mathbf{Q}_\ell \mathbf{Q}_\ell^T)\mathbf{b}$.  We note that a beneficial side effect of applying this preconditioner is that it deflates the zero eigenvalue associated with Neumann type boundary conditions, which would be required anyway before applying the contour integral method (since the method applies only to matrices whose eigenvalues all lie in the right half plane).

The inner level of preconditioning operates only at the level of the shifted linear system solves in \eqref{eq:contoursol}.  Care is required when constructing preconditioners for shifted systems, since the subspace is built using only the unshifted matrix $\mathbf{A}$.  One can show that certain polynomial preconditioners remain applicable to the family of shifted systems \cite{ahmad12preconditioned, farquhar16GPU}, and in this work we use the least squares polynomials obtained for the Jacobi weights with parameters $\mu = \frac{1}{2}$ and $\nu = -\frac{1}{2}$. The derivation of these polynomials and a table of examples can be found in Saad \cite[p.~400-402]{saad03iterative}.

\section{Numerical Results}
\label{sec:results4}
\subsection{One-dimensional Fisher problem}
\label{sec:nicoleres}
We begin this section by comparing the results of our approach outlined in the preceding sections to the results published in Cusimano~\cite{cusimano15fractional}.  In that work, a dense matrix representation of the variable-order two-sided Riemann-Liouville operator with reflecting boundary conditions was derived and utilised.  The computer code for generating the matrices for one-dimensional problems in this manner was generously provided by Cusimano, and so we are able to directly compare her approach with that proposed in this paper.

The test problem considered by Cusimano~\cite{cusimano15fractional} is a variable-order fractional reaction-diffusion problem
\begin{equation}
\frac{\partial u(x,t)}{\partial t} = - \left(-\nabla^2\right)^{\alpha(x)/2}u(x,t) + g(u)
\label{eq:test1}
\end{equation}
with a Fisher source term
\begin{equation*}
g(u) = u(1-u),
\end{equation*}
and piecewise-constant fractional order
\begin{equation*}
\alpha(x) = \left\{ \begin{aligned}
\alpha_1\,, & \ \ \  x \in \Omega_1 \\
\alpha_2\,, & \ \ \  x \in \Omega_2
\end{aligned}\right.
\end{equation*}
with $\alpha_1 = 1.5$, $\alpha_2 = 2$.  The domain of the problem is the closed interval $\Omega = [0,100]$ with the subdomains $\Omega_1 = [0,50]$ and $\Omega_2 = (50,100]$.  A step-like initial condition is used,
\begin{equation*}
u(x,0) = \left\{\begin{aligned}
&1 & x\leq 5\\
&e^{-10(x-5)} & x> 5.
\end{aligned}
\right.
\end{equation*}
To match the discretisation of Cusimano \cite{cusimano15fractional} we discretise using standard uniform finite differences with Neumann boundary conditions to generate the matrix representation $\mathbf{A}$ of the Laplacian operator.

\begin{figure}
	\centering
	\subfloat[$\alpha_1 = 1.5, \alpha_2 = 2$.]{\label{fig:1dres1}\includegraphics[width=0.49\textwidth]{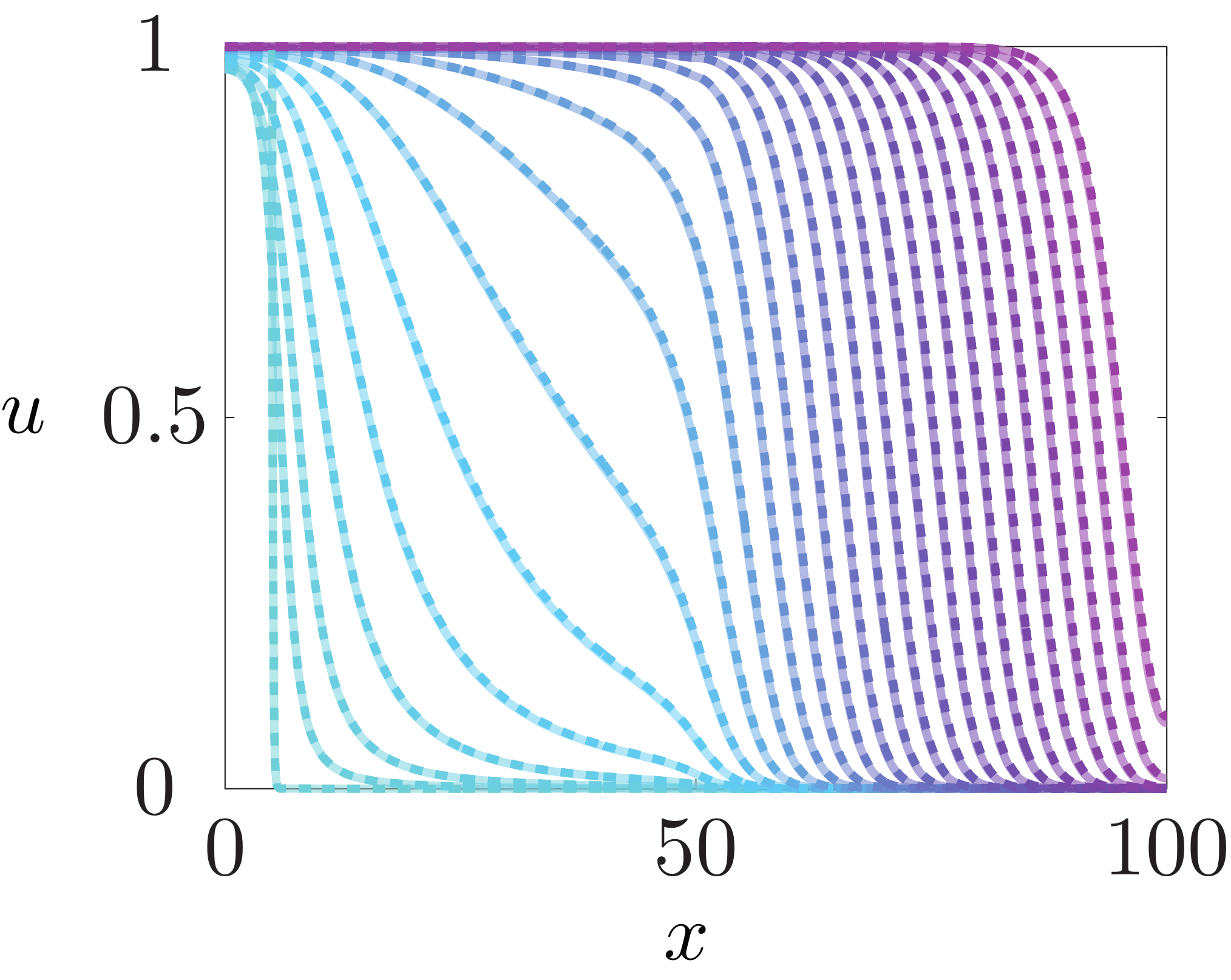}} \enskip
	\subfloat[$\alpha_1 = 2, \alpha_2 = 1.5$.]{\label{fig:1dres2}\includegraphics[width=0.49\textwidth]{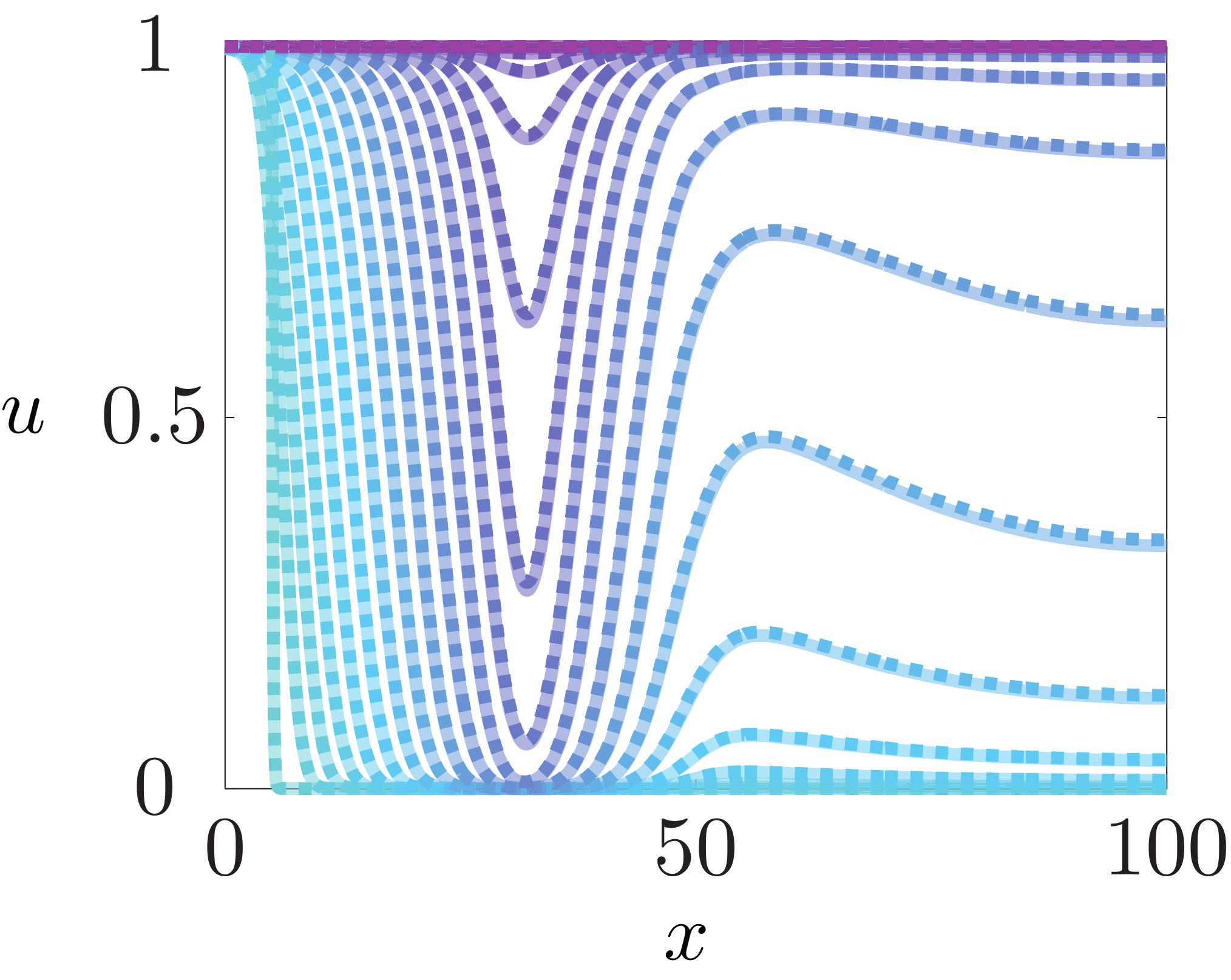}} \\
	\subfloat[$\alpha_1 = 1.5, \alpha_2 = 1.8$.]{\label{fig:1dres3}\includegraphics[width=0.49\textwidth]{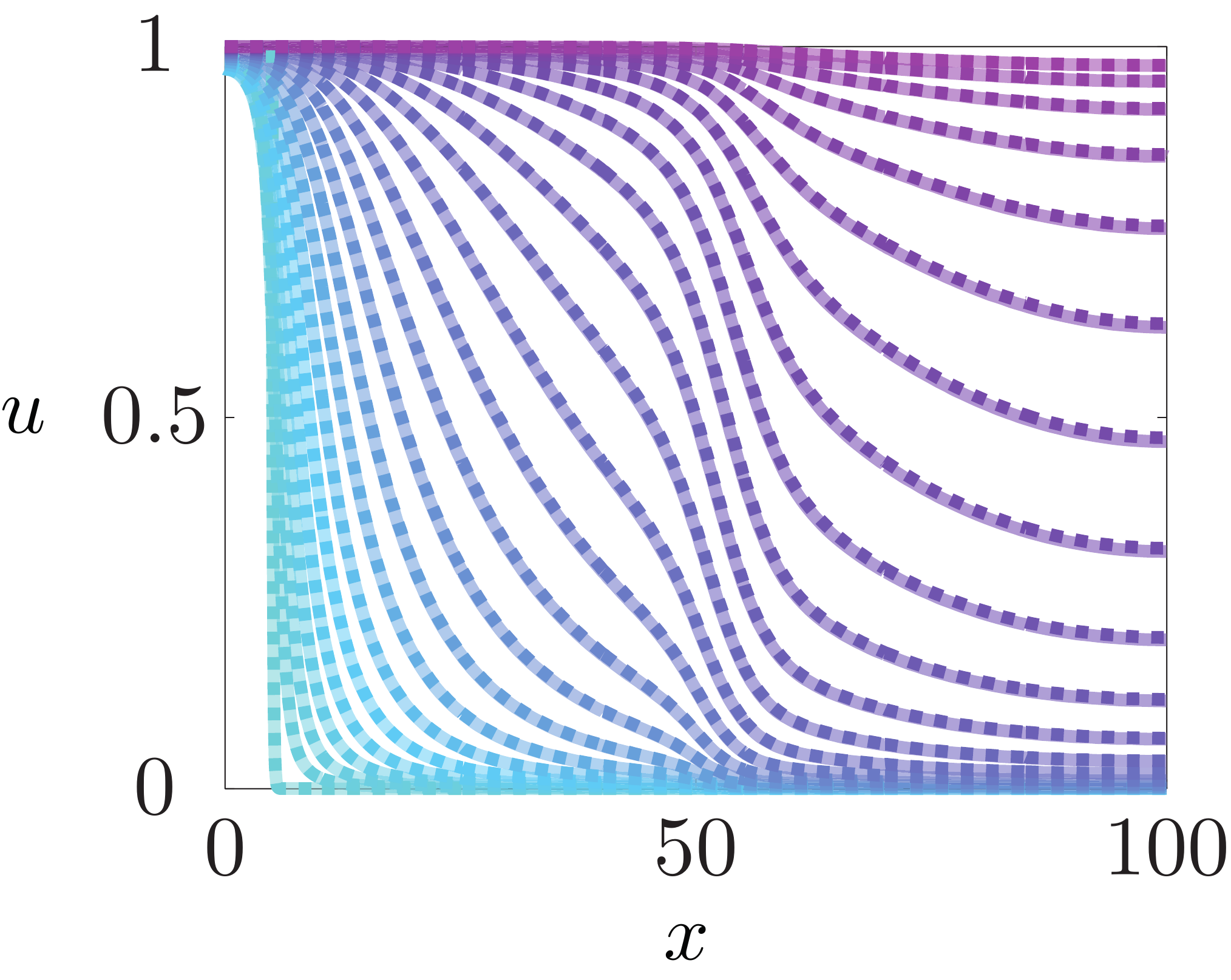}} \enskip
	\subfloat[$\alpha_1 = 1.8, \alpha_2 = 1.5$.]{\label{fig:1dres4}\includegraphics[width=0.49\textwidth]{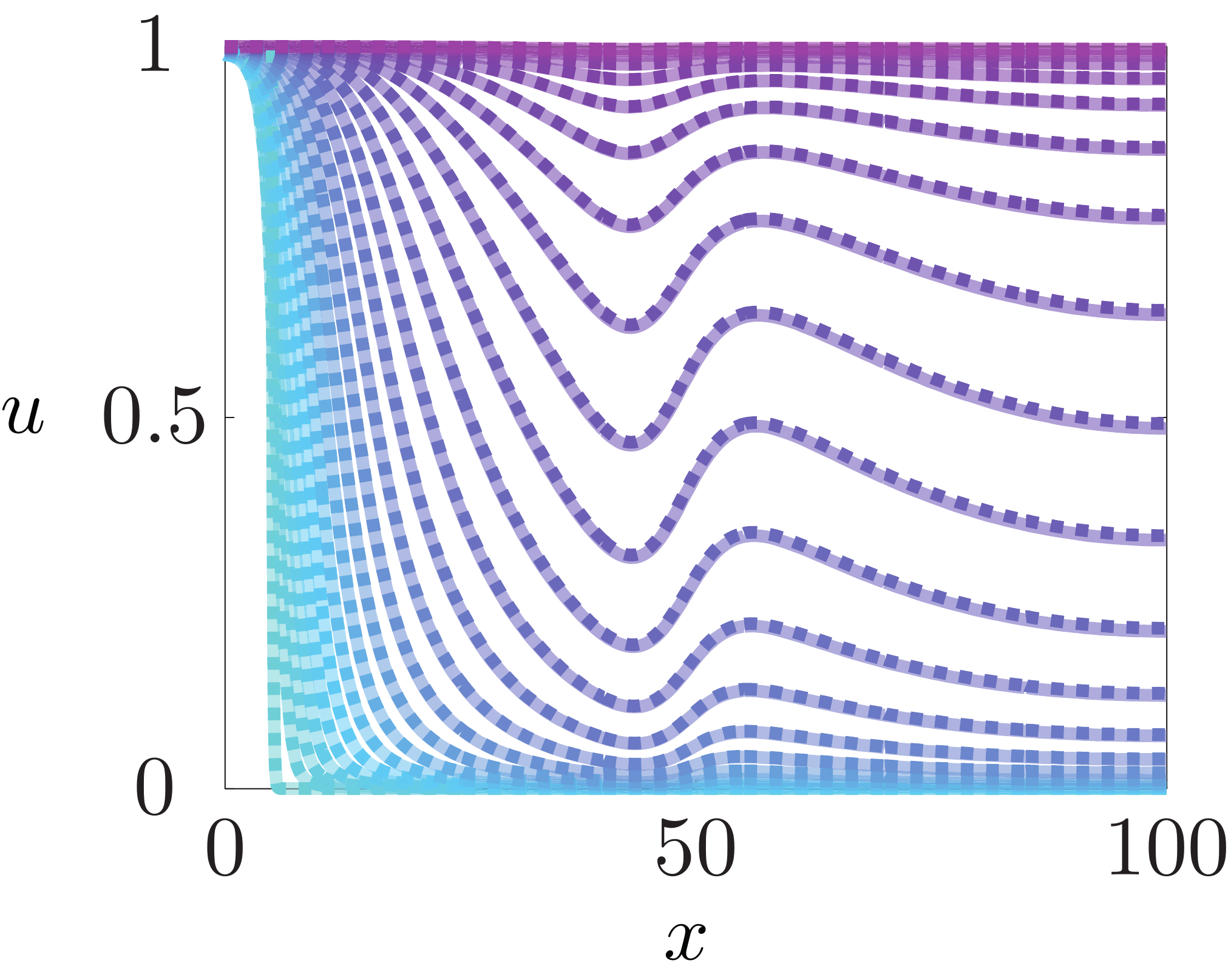}}
	\caption{Comparison of the solutions to Equation~\eqref{eq:test1} using the method of Cusimano~\cite{cusimano15fractional} (solid lines) against the variable-order fractional Laplacian (dashed lines) for different combinations of $\alpha_1$ and $\alpha_2$.}
\end{figure}

The results for this problem can be seen in Figure~\ref{fig:1dres1}.  Solid curves are the solutions using the method of Cusimano, and dashed curves indicate our present solutions, with no disagreement between solutions observed.
They exhibit the characteristic exponential-type acceleration of the solution to the $u = 1$ state in the left half of the domain where $\alpha_1 = 1.5$, before transitioning into travelling wave behaviour in the right half of the domain where $\alpha_2 = 2$.

Using the provided code, we also considered a variety of different choices for $\alpha(x)$ than those published in Cusimano~\cite{cusimano15fractional}. In Figure~\ref{fig:1dres2}, we present the results for $\alpha_1 = 2$ and $\alpha_2 = 1.5$. We observe that in the left half of the domain a travelling wave solution develops whose profile is ultimately affected by the smaller value of $\alpha$ in the right half of the domain, forcing the solution to transition to a regime that approaches the steady state exponentially fast \cite{cusimano15fractional}.

In Figure~\ref{fig:1dres3}, we present the results for $\alpha_1 = 1.5$ and $\alpha_2 = 1.8$. We see the acceleration towards steady-state is more rapid in the left half of the domain, associated with $\alpha = 1.5$ than in the right half of the domain with $\alpha = 1.8$.  The results for the final combination of parameters, that is $\alpha_1 = 1.8$, $\alpha_2 = 1.5$, are given in Figure~\ref{fig:1dres4}, where now the acceleration towards steady-state is faster in the right half of the domain, consistent with the smaller value of $\alpha$ in that region.

In all figures we have compared the results from our method, which uses sparse matrices and matrix function vector product techniques, with those using Cusimano's method, which uses dense matrices, and we observe complete agreement in all cases.

Simulations were run on a Dell desktop machine with Windows 7, Intel(R) Core(TM) i7-4770 CPU \@ 3.40GHz Processor, 16.0 GB RAM using MATLAB R2016a.
In Table~\ref{tab:compnic} we present comparisons of a number of important characteristics concerning the efficiency of each method, including the memory storage and runtime.  These results show that the iteration time for the method of Cusimano is much faster than our method utilising the variable-order fractional Laplacian (VOFL) for small matrix sizes. As the problem size is increased, this discrepancy in runtimes diminishes, and by dimension $N = 4001$ the runtimes are within a factor of two of each other.  In fact, for the special case where $\alpha_1 = 2$, one can recover even this factor of two by recognising that the matrix function $f_b(\mathbf{A}) = \left(\mathbf{I} + D\delta t \mathbf{A}^{\alpha_{1}/2} \right)^{-1}$ in \eqref{eq:BEfp} is now simply a linear system solve, which does not require the full machinery of matrix function vector product evaluation.  (In the case that $\alpha_2 = 2$ one can simply swap the roles of $\alpha_1$ and $\alpha_2$ to exploit this same efficiency gain.)

A further requirement for the method of Cusimano is the actual formation of the dense coefficient matrix.  We list the runtime for this step separately in Table~\ref{tab:compnic} since the efficient construction of this matrix was not the focus of her thesis.  However, in the cases tested, this step actually accounted for the vast majority of the total runtime.  No such computational effort is required for our proposed method, since we deal only with the sparse matrix representation of the standard Laplacian.

Another important factor to consider is the memory storage required to solve the problem.  The method of Cusimano requires much higher memory costs because of the requirement to store the dense coefficient matrix. For the one-dimensional problems presented, the mesh sizes are not refined to an extent that this becomes a problem, however memory storage constraints would emerge with further refinement or extension to higher dimensions.  In any case, the method of Cusimano does not easily extend to higher dimensions and especially to irregular domains, owing to the complex (but ingenious) method used to incorporate reflecting boundary conditions.  Our approach by contrast readily generalises to arbitrary dimensions and geometries, and the sparse matrix implementation avoids the memory constraints that otherwise limit dense matrix-based approaches.

\begin{table}
	\centering
	\footnotesize
	\begin{tabular}{>{\centering\arraybackslash}p{0.048\textwidth}	>{\centering\arraybackslash}p{0.09\textwidth}	>{\centering\arraybackslash}p{0.106\textwidth}	>{\centering\arraybackslash}p{0.088\textwidth}	>{\centering\arraybackslash}p{0.024\textwidth}	>{\centering\arraybackslash}p{0.031\textwidth}	>{\centering\arraybackslash}p{0.031\textwidth}	>{\centering\arraybackslash}p{0.093\textwidth}	>{\centering\arraybackslash}p{0.091\textwidth}	>{\centering\arraybackslash}p{0.04\textwidth}	}
		\toprule										
		$N$&	Number of Steps&	Operator&	 Vectors in \newline Memory&	$T$&	$\alpha_1$&	$\alpha_2$&	Iteration Time (s)&	Matrix Gen Time (s)&	Av. FP	\\
		\midrule										
		\multirow{8}{*}{$1001$}&	\multirow{8}{*}{$3000$}&	\multirow{4}{*}{Cusimano}&	\multirow{4}{*}{$1001$}&	\multirow{2}{*}{$15$}&	\mdef 1.5&	\mdef 1.8&	\mdef 4.8&	\mdef 281.6&	\mdef 	\\
		&	&	&	&	&	1.8&	1.5&	4.7&	281.6&		\\
		&	&	&	&	\multirow{2}{*}{$30$}&	\mdef 1.5&	\mdef 2&	\mdef 4.8&	\mdef 240.6&	\mdef 	\\
		&	&	&	&	&	 2&	1.5&	4.8&	241.9&		\\
		\cmidrule{3-10}										
		&	&	\multirow{4}{*}{VOFL}&	\multirow{4}{*}{$156$}&	\multirow{2}{*}{$15$}&	\mdef 1.5&	\mdef 1.8&	\mdef 54.0&	\mdef &	\mdef 2.0	\\
		&	&	&	&	&	 1.8&	1.5&	51.2&	&	1.8	\\
		&	&	&	&	\multirow{2}{*}{$30$}&	\mdef 1.5&	 \mdef 2&	\mdef 26.9&	\mdef &	\mdef 2.1	\\
		&	&	&	&	&	 2&	1.5&	21.6&	&	1.7	\\
		\midrule										
		\multirow{8}{*}{$2001$}&	\multirow{8}{*}{$12000$}&	\multirow{4}{*}{Cusimano}&	\multirow{4}{*}{$2001$}&	\multirow{2}{*}{$15$}&	\mdef 1.5&	\mdef 1.8&	\mdef 76.5&	\mdef 2273.6&	\mdef 	\\
		&	&	&	&	&	1.8&	1.5&	72.0&	2267.7&		\\
		&	&	&	&	\multirow{2}{*}{$30$}&	\mdef 1.5&	\mdef 2&	\mdef 75.0&	\mdef 1943.7&	\mdef 	\\
		&	&	&	&	&	2&	1.5&	76.7&	1951.2&		\\
		\cmidrule{3-10}										
		&	&	\multirow{4}{*}{VOFL}&	\multirow{4}{*}{$221$}&	\multirow{2}{*}{$15$}&	\mdef 1.5&	\mdef 1.8&	\mdef 218.0&	\mdef &	\mdef 1.0	\\
		&	&	&	&	&	 1.8&	1.5&	218.6&	&	1.0	\\
		&	&	&	&	\multirow{2}{*}{$30$}&	\mdef 1.5& \mdef 2&	\mdef 112.0&	\mdef &	\mdef 1.1	\\
		&	&	&	&	&	 2&	1.5&	113.0&	&	1.0	\\
		\midrule										
		\multirow{8}{*}{$4001$}&	\multirow{8}{*}{$48000$}&	\multirow{4}{*}{Cusimano}&	\multirow{4}{*}{$4001$}&	\multirow{2}{*}{$15$}&	\mdef 1.5&	\mdef 1.8&	\mdef 1041.0&	\mdef 17962.4&	\mdef 	\\
		&	&	&	&	&	1.8&	1.5&	1042.6&	17963.4&		\\
		&	&	&	&	\multirow{2}{*}{$30$}&	\mdef 1.5&	\mdef 2&	\mdef 1075.1&	\mdef 15327.6&	\mdef 	\\
		&	&	&	&	&	2&	1.5&	1049.1&	15324.0&		\\
		\cmidrule{3-10}										
		&	&	\multirow{4}{*}{VOFL}&	\multirow{4}{*}{$312$}&	\multirow{2}{*}{$15$}&	\mdef 1.5&	\mdef 1.8&	\mdef 2253.4&	\mdef &	\mdef 1.0	\\
		&	&	&	&	&	 1.8&	1.5&	2302.4&	&	1.0	\\
		&	&	&	&	\multirow{2}{*}{$30$}&	\mdef 1.5& \mdef 2&	\mdef 1090.9&	\mdef &	\mdef 1.0	\\
		&	&	&	&	&	 2&	1.5&	1115.5&	&	1.0	\\
		\bottomrule
	\end{tabular}
	\caption[Table of timing and storage comparisons between the method of Cusimano~\cite{cusimano15fractional} and the variable-order fractional Laplacian.]{Table of timing and storage comparisons between the operator of Cusimano~\cite{cusimano15fractional} and the variable-order fractional Laplacian (VOFL).  $N$ is the number of nodes in the domain, Vectors in Memory is the number of vectors $\mathbb{R}^{N\times 1}$ that are required to be stored, Matrix generation time is the time in seconds to generate the matrix for the operator of Cusimano \cite{cusimano15fractional} and Av. FP is the average number of fixed point iterations required for the scheme involving the variable-order fractional Laplacian. }
	\label{tab:compnic}
\end{table}

\subsection{One-dimensional Beeler-Reuter model}
\label{sec:1dBRres}
\begin{figure}
	\centering
	{\includegraphics[width = 0.85\textwidth]{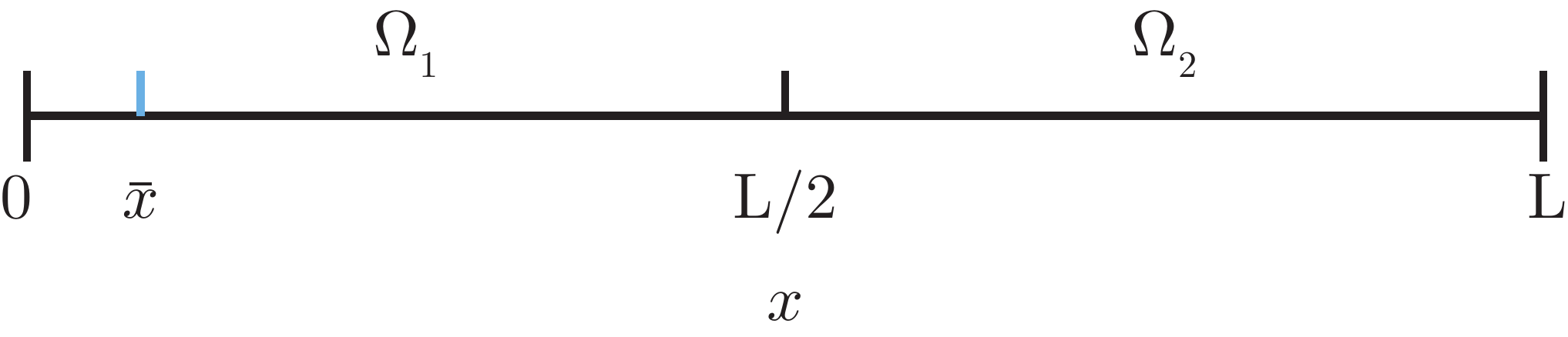}}
	\caption[Solution domain considered in one dimension for the coupled monodomain, ionic current models.]{Solution domain considered in one dimension for the coupled monodomain, ionic current models. This domain represents a cable of tissue, with a stimulus applied to the left end of the cable, between $0$ and the blue line $\bar{x}$.}
	\label{fig:1ddomain5}
\end{figure}

Moving now to cardiac models, in one dimension we further consider the variable-order fractional monodomain equation \eqref{eq:varfracmono},
coupled with the Beeler-Reuter model for ionic currents \eqref{eq:BRode}.  This one-dimensional test problem simulates the propagation of an electrical impulse along a cable of tissue in the heart.  The domain that we consider is $\Omega = [0,L]$, as shown in Figure~\ref{fig:1ddomain5}, split in the middle such that for the variable-order operator we have
\begin{equation*}
\alpha(x) = \left\{ \begin{aligned}
\alpha_1 & \ \ \ x \in \Omega_1 = [0,L/2]\\
\alpha_2 & \ \ \ x \in \Omega_2 = (L/2,L].
\end{aligned}\right.
\end{equation*}
These simulations continue to use a uniform finite difference spatial discretisation of the Laplacian operator. The stimulus current is applied to an interval $\Omega_{stim}=[0,\bar{x}]$, between the left end of the domain and $\bar{x}$ as marked in Figure~\ref{fig:1ddomain5}.

We have chosen monodomain parameters similar to those of Cusimano~\cite{cusimano15fractional}, $C_m = 1 \mu \textrm{F}\cdot \textrm{cm}^{-2}$, $\chi = 2000 \textrm{cm}^{-1}$ and  $D = 1 \textrm{mS}\cdot\textrm{cm}^{-1}$. We simulate the problem to an end time of 1200\textrm{ms} taking time steps of $\delta t = 0.25$. We set $L = 10$ and use a node spacing of $\Delta x = 0.01$. For the value of the applied stimulus current we follow the example of Cusimano et al.~\cite{cusimano15order} and set the value of the applied stimulus current to be twice the diastolic threshold in the standard diffusion case, which is defined to be the minimum current that produces successful propagation of the electrical impulse. To allow the gating variables and the electrical potential to reach a rest state, the stimulus current is applied after $T_{stim} = 10 \textrm{ms}$ for a continuous interval of $\delta t_{stim} = 5\textrm{ms}$ on the spatial interval $\Omega_{stim} = [0,\bar{x}] = [0,0.25\textrm{cm}]$. For the parameter values outlined here, the applied stimulus current is given by $I_{stim} = 12\chi$. The initial conditions for this problem are set such that the transmembrane potential has a value of $v = -85\textrm{mV}$ at every node in the domain and the gating variables have the initial values $m = 0$, $h=1$, $j= 1$, $d=0$, $f=1$, $x = 0$ and  $c = 1$ everywhere in the domain.

\begin{figure}
	\centering
	\subfloat[$\alpha = 2$]{\label{fig:BRtw1dstandard}\includegraphics[width = 0.48\textwidth]{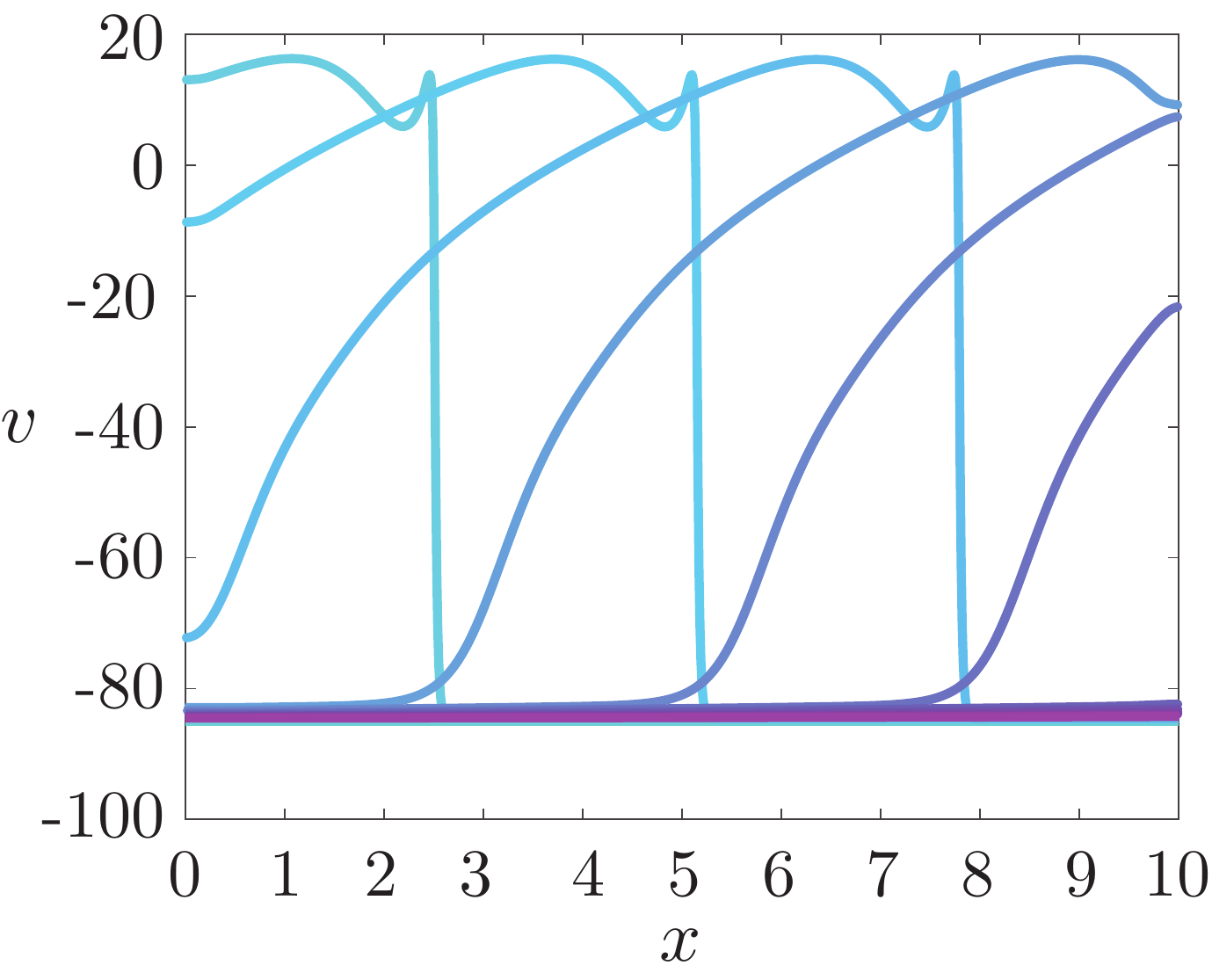}}\enskip
	\subfloat[$\alpha = 1.5$]{\label{fig:BRtw1dfrac}\includegraphics[width = 0.48\textwidth]{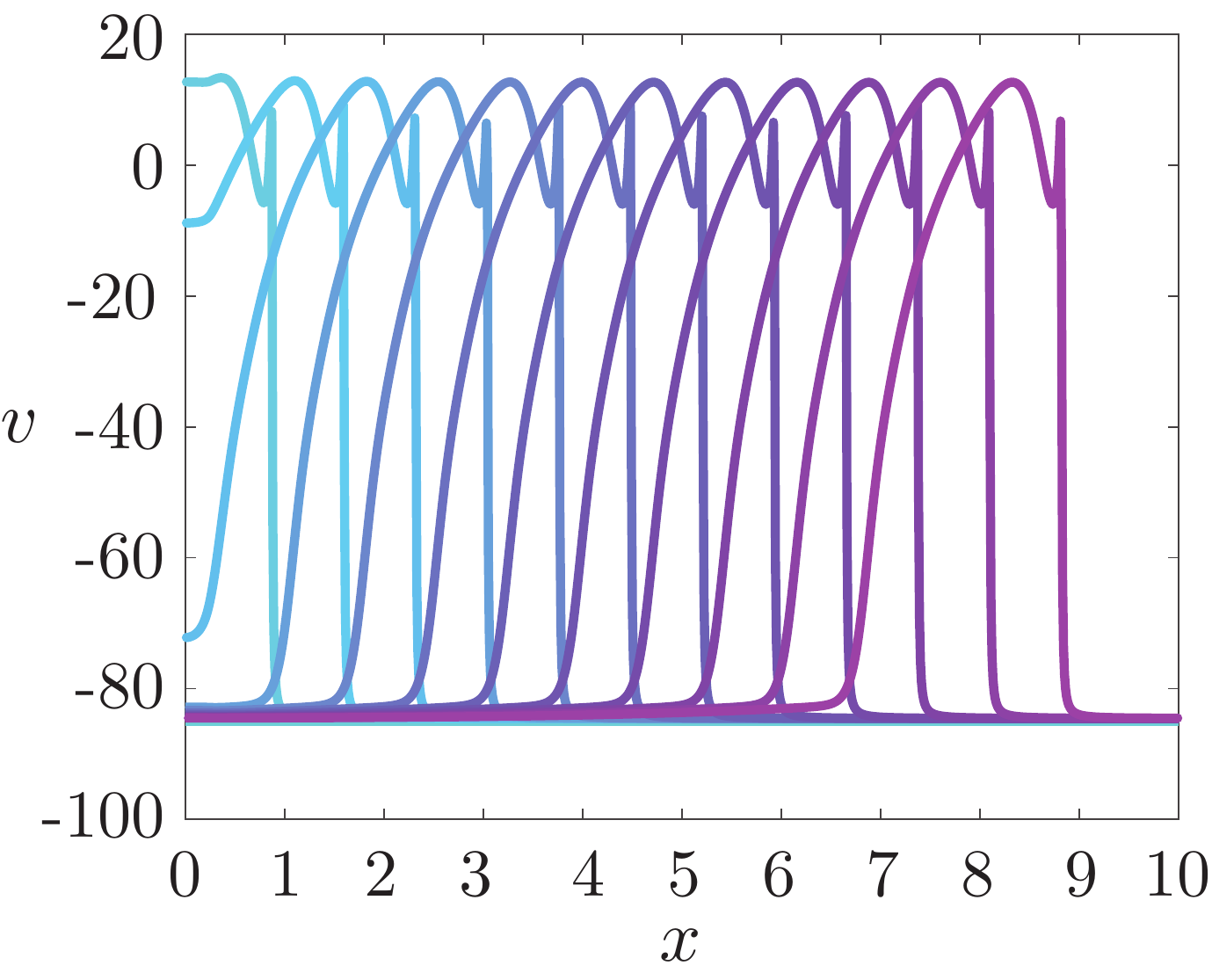}}\\
	\caption{Solution profile of the Beeler-Reuter monodomain model in one-dimension with (a) standard diffusion ($\nabla^2$) and (b) fractional diffusion ($-(-\nabla^2)^{\alpha/2}$) as a function of $x$.}
	\label{fig:BRtravellingwaves1d}
\end{figure}

We first computed the solution using both standard diffusion ($\alpha_1 = \alpha_2 = 2$) and fixed-order fractional diffusion ($\alpha_1 = \alpha_2 = 1.5$). The waves of these problems are shown in Figure~\ref{fig:BRtravellingwaves1d}, with the wave profile plotted for every $100 \textrm{ms}$. In these results we can see that, unlike the Fisher problem, the Beeler-Reuter model produces travelling wave type solutions for both the standard diffusion case (Figure~\ref{fig:BRtw1dstandard}) and the fractional diffusion case (Figure~\ref{fig:BRtw1dfrac}). We observe that the speed of the wave in the standard diffusion case is faster than in the fractional diffusion case, as is clear by the number of waves present in each plot. However the general shape of the wave profile is the same for both cases, consisting of a slow increase of the wave to a plateau before a smooth drop that creates a notch when the wave spikes, before finally falling rapidly to steady state to create an almost vertical front.

\begin{figure}
	\centering
	\subfloat[][$\alpha(x) =\left\{\begin{aligned}& 1.5, \ \  &x\leq L/2\\ & 2, \ \  &x>L/2
	\end{aligned}\right.$]{\label{fig:BRtwvolh}\includegraphics[width = 0.48\textwidth]{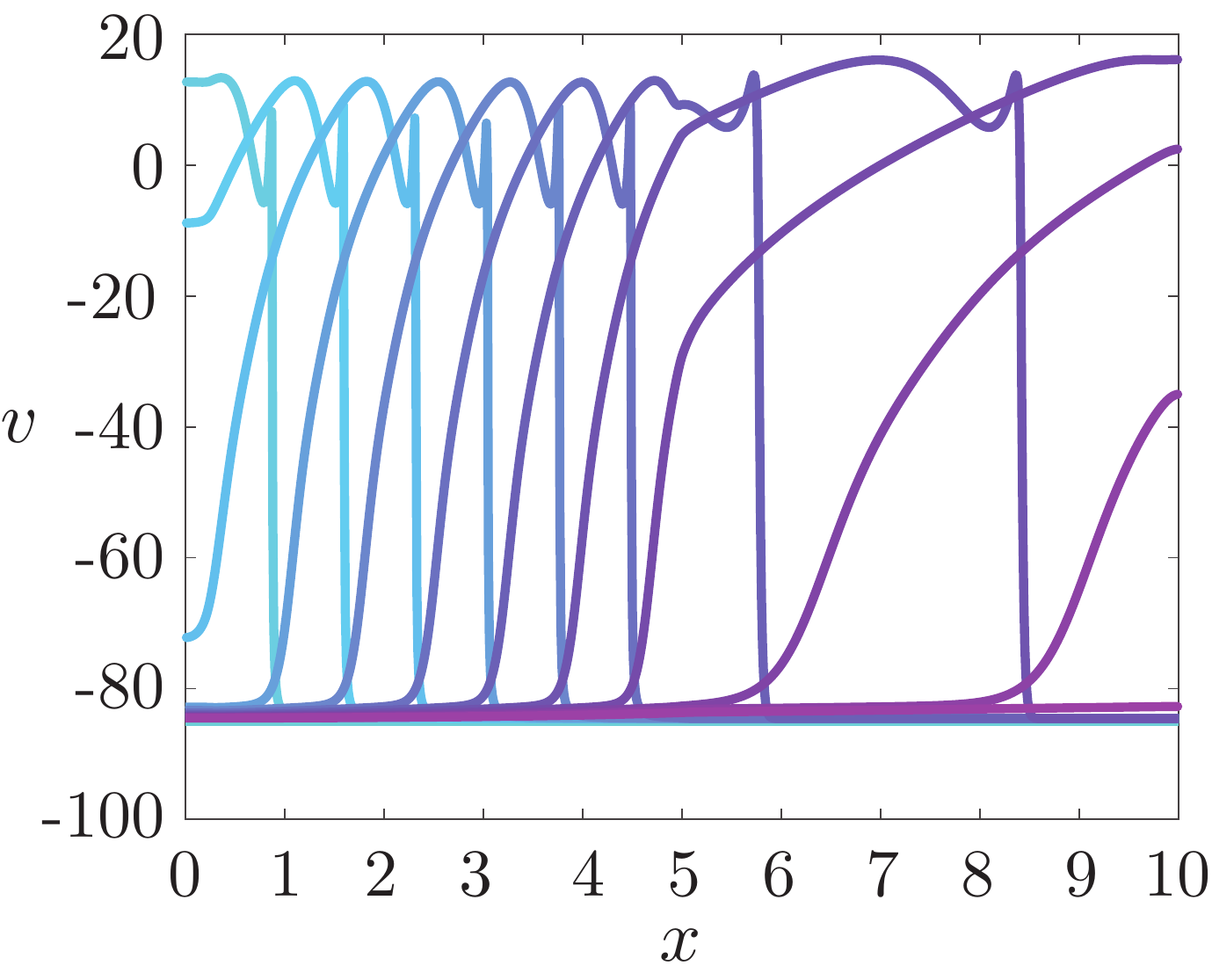}}\enskip
	\subfloat[][$\alpha(x) =\left\{\begin{aligned}& 2, \ \  &x\leq L/2\\ & 1.5, \ \  &x>L/2
	\end{aligned}\right.$]{\label{fig:BRtwvohl}\includegraphics[width = 0.48\textwidth]{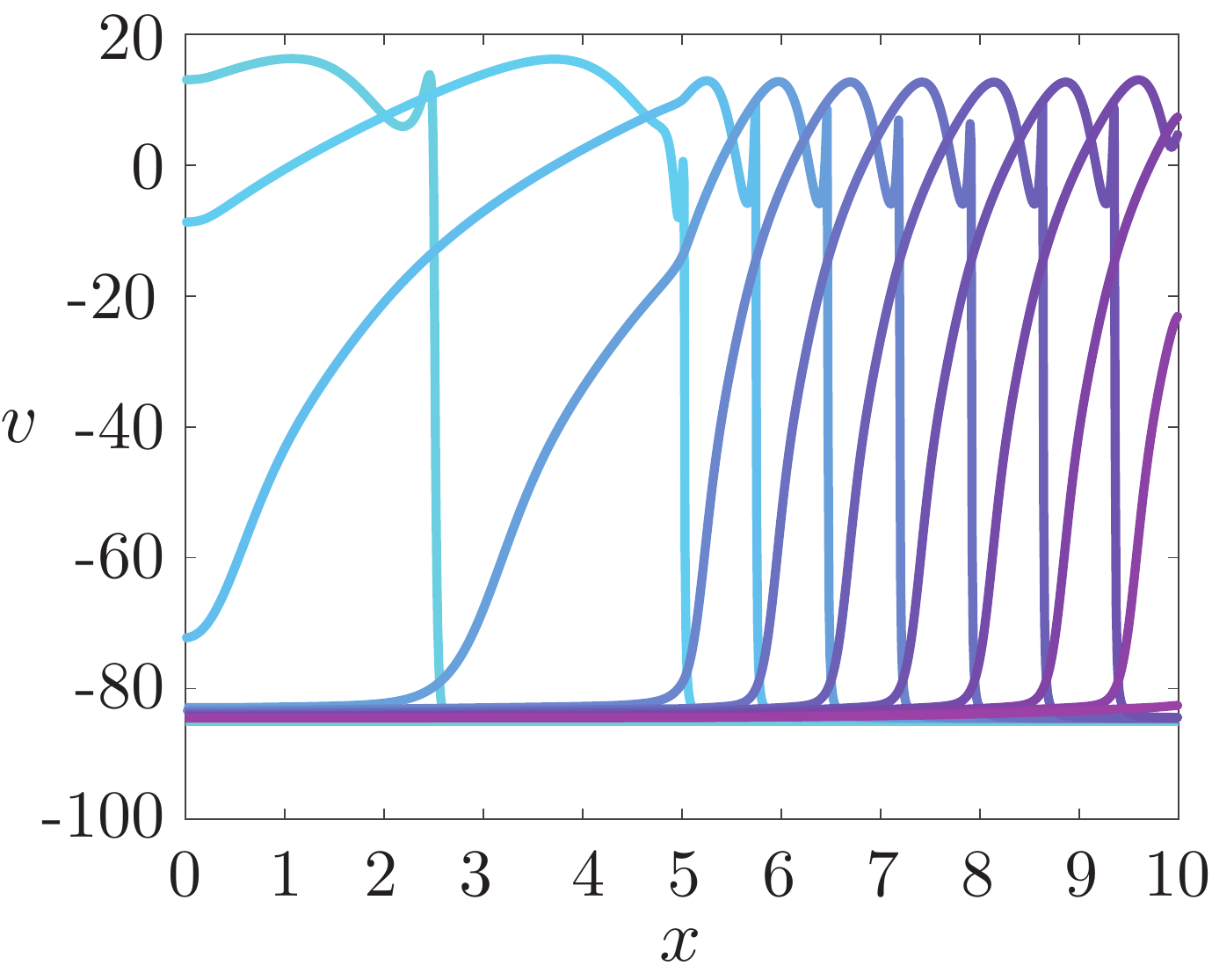}}
	\caption{One-dimensional solutions to the Beeler-Reuter monodomain model with variable-order fractional diffusion ($-(-\nabla^2)^{\alpha(x)/2}$).}
	\label{fig:BR1dvartw}
\end{figure}

Results that compare the effect of different fractional indices on the wave profile generated by the fixed-order fractional monodomain model coupled with the Beeler-Reuter ionic current model have been previously presented in Cusimano~\cite{cusimano15fractional} and Cusimano et al.~\cite{cusimano15order}. Our contribution is to extend the model to incorporate the variable-order fractional Laplacian operator, which accounts for the differences in the tissue structure in two regions.

In Figure~\ref{fig:BR1dvartw} we present the results that incorporate the variable-order fractional Laplacian operator. Figure~\ref{fig:BRtwvolh} exhibits the result with $\alpha_1 = 1.5$ and $\alpha_2 = 2$ . In this figure we see that the wave resembles a travelling wave that changes form as it passes through the boundary between the regions. In the left half of the domain, $x \in \Omega_1$, the wave resembles the shape and speed of the travelling waves in Figure~\ref{fig:BRtw1dfrac}. That is, the wave for $x \in\Omega_1$ resembles the fixed-order fractional problem with $\alpha = \alpha_1 = 1.5$. In the right half of the domain, the shape and speed of the wave instead resembles the travelling waves in Figure~\ref{fig:BRtw1dstandard}. That is, for $x\in \Omega_2$, the wave resembles the standard order problem, $\alpha = \alpha_2 = 2$.  In Figure~\ref{fig:BRtwvohl} we reverse the values of $\alpha_1$ and $\alpha_2$, and the conclusions hold similarly: the results are consistent with those of Figure~\ref{fig:BRtw1dfrac} for the corresponding values of $\alpha$, with a transition in the middle of the domain.

\subsection{Three-Dimensional Problem - Modelling Ischaemic re-entry with Fractional Diffusion}
\label{sec:3dBRres}
We now conduct full three-dimensional simulations of the identical model consisting of the variable-order fractional monodomain equation \eqref{eq:varfracmono} coupled with the Beeler-Reuter model \eqref{eq:BRode}.  We use the mesh of a rabbit heart from the CHASTE package~\cite{ChasteGRM2013,ChasteJPF2009}: specifically, the UCSD unstructured tetrahedral mesh consisting of $63885$ nodes and $322267$ elements that can be downloaded from \url{https://chaste.cs.ox.ac.uk/}.  The dimensions of the mesh are approximately $3.31 \textrm{cm} \times 2.92 \textrm{cm} \times 2.75 \textrm{cm}$.

The spatial discretisation for the Laplacian operator is generated using a vertex-centred finite volume method. In this case, the discretisation produces a diagonal mass matrix $\mathbf{M}$ and a stiffness matrix $\mathbf{K}$ such that the matrix representation of the Laplacian is $\mathbf{A} = \mathbf{M}^{-1} \mathbf{K}$.  Working with such a matrix in the framework presented requires one additional step~\cite{farquhar16GPU,yang10novel}: define $\tilde{\mathbf{A}} = \mathbf{M}^{-1/2} \mathbf{K} \mathbf{M}^{-1/2}$ to recover the required symmetric matrix $\tilde{\mathbf{A}}$, and compute $f(\mathbf{A}) \mathbf{b}$ by
\begin{equation*}
f(\mathbf{A}) \mathbf{b} = \mathbf{M}^{-1/2} f(\tilde{\mathbf{A}})\tilde{\mathbf{b}}, \ \ \ \ \tilde{\mathbf{b}} = \mathbf{M}^{-1/2}\mathbf{b}.
\end{equation*}

\begin{figure}
	\centering
	\includegraphics[trim=50 0 25 0,clip,width = 0.7\textwidth]{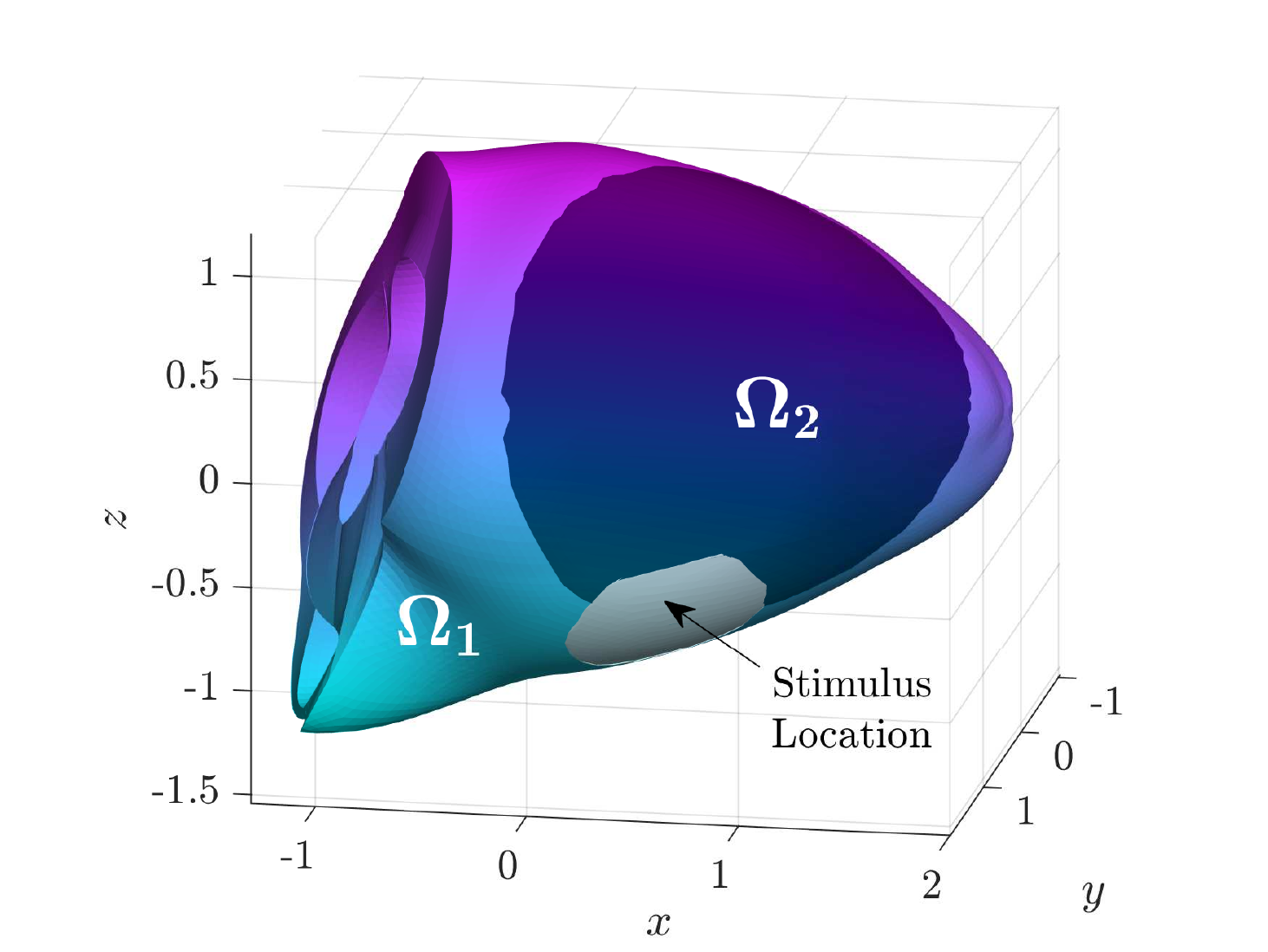}
	\caption[The solution domain considered in three dimensions for the coupled monodomain, ionic current model with re-entry]{The solution domain considered in three dimensions for the coupled monodomain, ionic current model with re-entry This domain is a mesh of a rabbit heart. The stimulus current is applied to the lighter region of the image, as marked. The region in darker regions represents ischaemic tissue as marked by $\Omega_2$.}
	\label{fig:heartmesh2}
\end{figure}

For the stimulus location we choose a partial sphere with radius 0.5 centred at the node $(0.3513, 0.0707, -1.0772)$. We model an ischaemic region that is a partial sphere centred at the node $(0.6079, -0.7698, -0.1947)$ with radius 1.25 excluding the tissue in between the ventricles,
\[\begin{aligned}
& &\Omega_2 = \Omega \cap \left\{ (x,y,z) | \sqrt{(x - 1.0352)^2 + (y+0.6256)^2 + (z - 0.248)^2}\leq 1.25 \right\}& & \\
& &\setminus \left\{ x<1.3, y>0.095, x>-0.3 \right\}& \,.
\end{aligned}\]
The remainder of the domain, $\Omega_1 = \Omega\setminus\Omega_2$, is composed of healthy tissue. A visual depiction of this domain can be seen in Figure~\ref{fig:heartmesh2}, where $\Omega_1$ is the region of healthy tissue and $\Omega_2$ is the ischaemic region.

We have used model parameters as for the one-dimensional simulation: $C_m~=~1~\mu~\textrm{F}\cdot \textrm{cm}^{-2}$, $\chi = 2000 \textrm{cm}^{-1}$ and  $D = 2 \textrm{mS}\cdot\textrm{cm}^{-1}$. The simulation is run to an end time of $1500\textrm{ms}$, taking time steps of $\delta t = 0.25\textrm{ms}$. We apply a repeated stimulus current (every $325\textrm{ms}$) at $T_{stim} = \{10, 335, 660, 985, 1310\}$ for a continuous period of $\delta t_{stim} = 5\textrm{ms}$ each time. {The choice of time between each stimulus is an important component to observing reentrant behaviour, too short a time between stimulus and the ``healthy'' region will not have recovered to be reactivated, while too long a time between stimulus and the ``damaged'' region will have recovered enough to be reactivated with the new stimulus.} Results are shown for a stimulus current applied to the region marked ``Stimulus location'' in Figure~\ref{fig:heartmesh2}. As with the previous test problem, the value of the stimulus current was chosen to be twice the diastolic threshold, which was $I_{stim} = 14\chi$. The initial conditions for this problem are that the transmembrane potential is $v = -85\textrm{mV}$ at every node in the domain, and the gating variables and calcium have the values $m = 0$, $h=1$, $j= 1$, $d=0$, $f=1$, $x = 0$, $c = 1\textrm{mol}\cdot \textrm{l}^{-1}$ at every node in the domain. In these results we consider the case where the fractional index is given by
\[\alpha(x) = \left\{ \begin{aligned} 2 & \ \ \ \ x \in \Omega_1 \\ 1.7 & \ \ \ \ x \in \Omega_2 \end{aligned}\right..\]
\begin{figure}
	\centering
	\subfloat[$t=20$]{\includegraphics[trim=140 40 10 50, clip, width=0.49\textwidth]{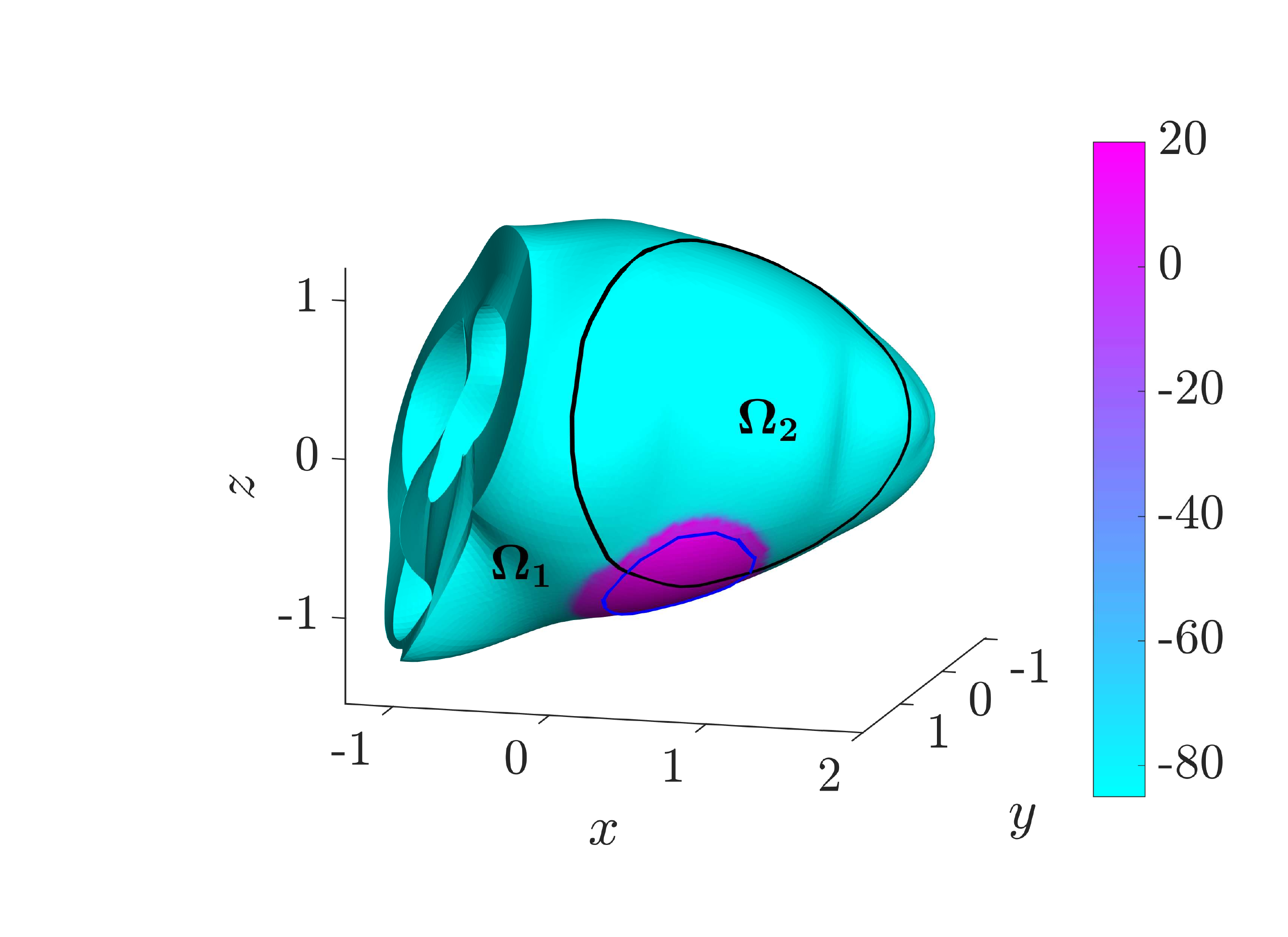}} \enskip
	\subfloat[$t=50$]{\includegraphics[trim=140 40 10 50, clip,width=0.49\textwidth]{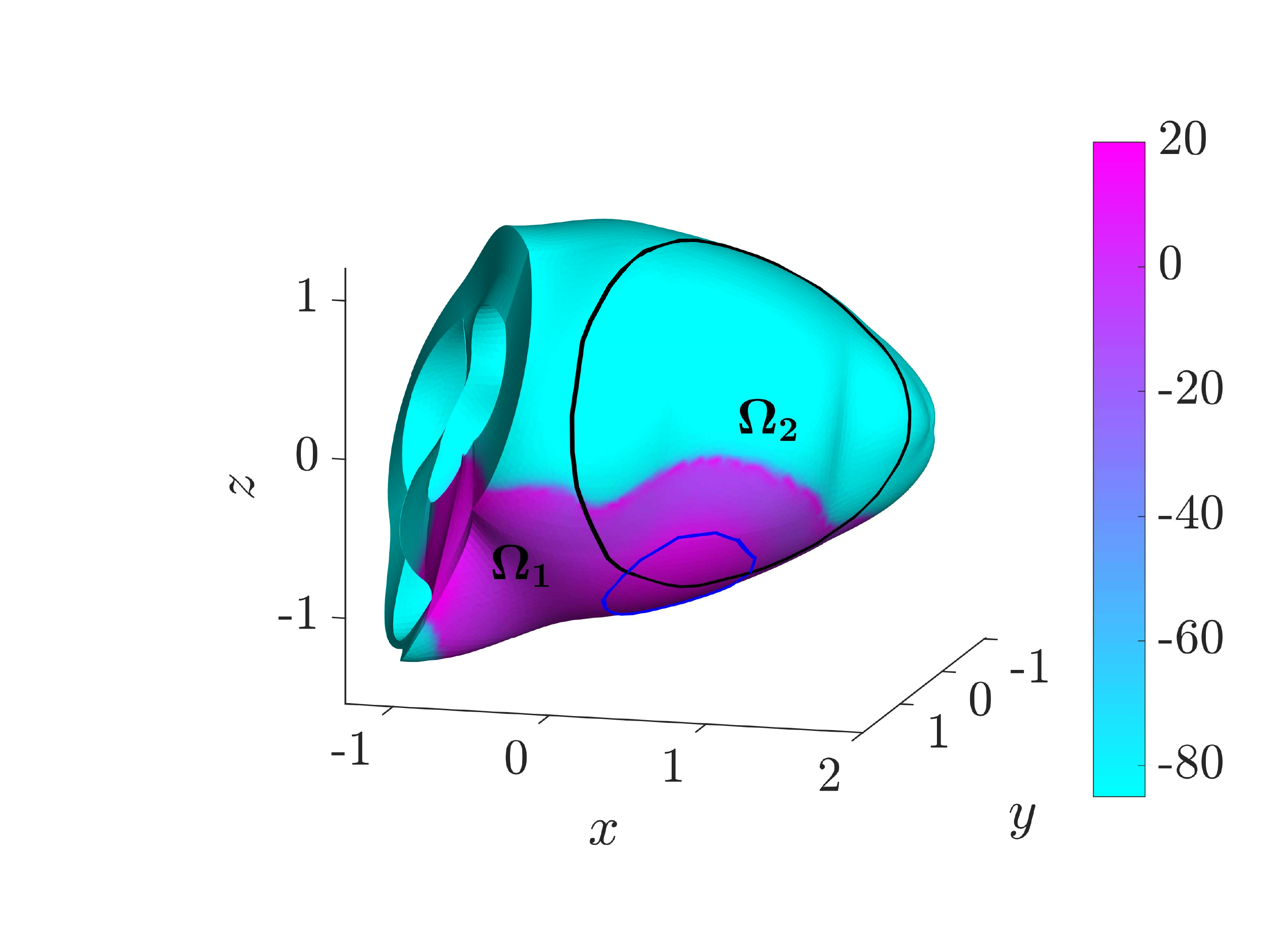}}
	\\
	\subfloat[$t=80$]{\includegraphics[trim=140 40 10 50, clip,width=0.49\textwidth]{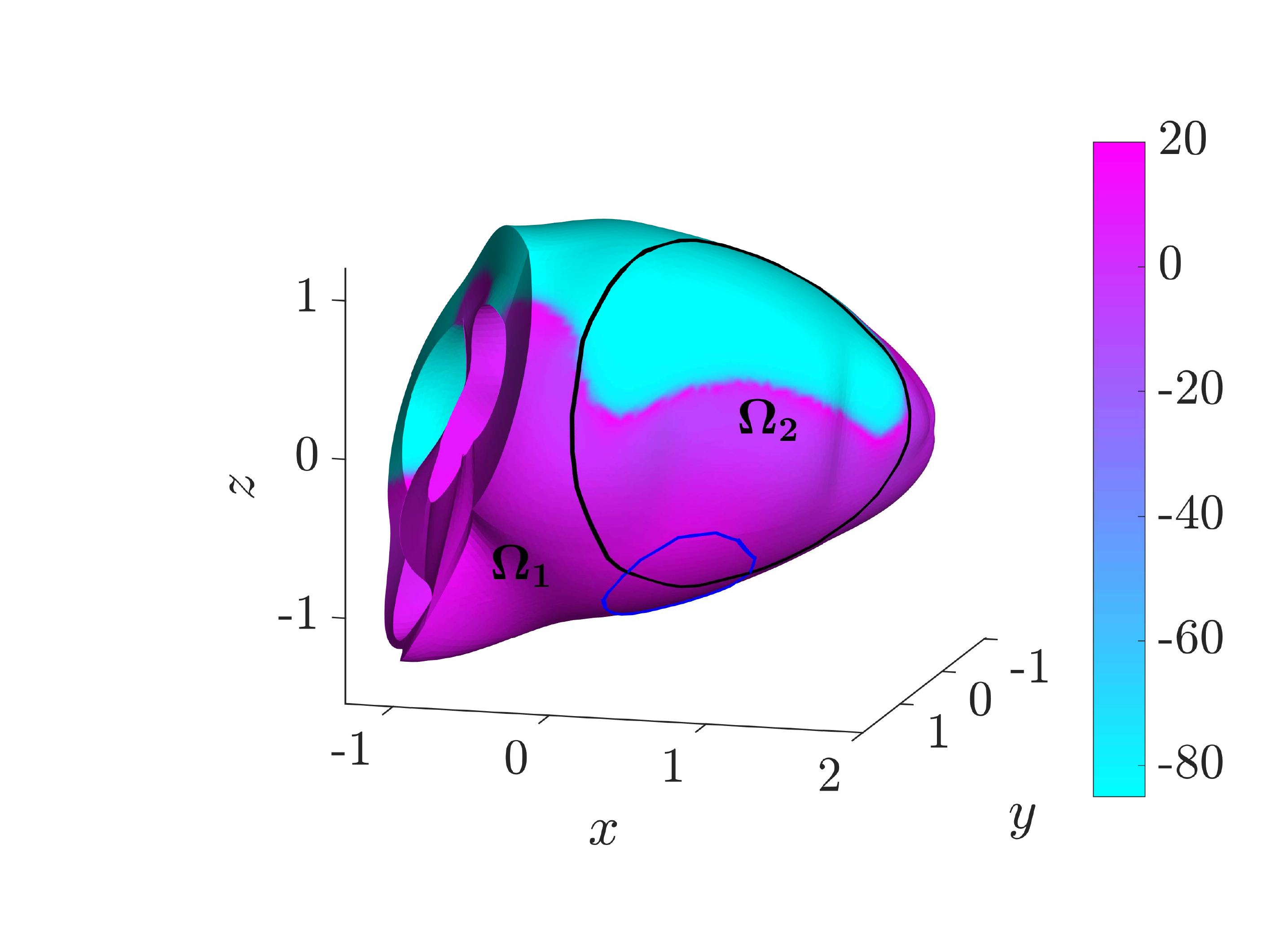}} \enskip
	\subfloat[$t=110$]{\includegraphics[trim=140 40 10 50, clip,width=0.49\textwidth]{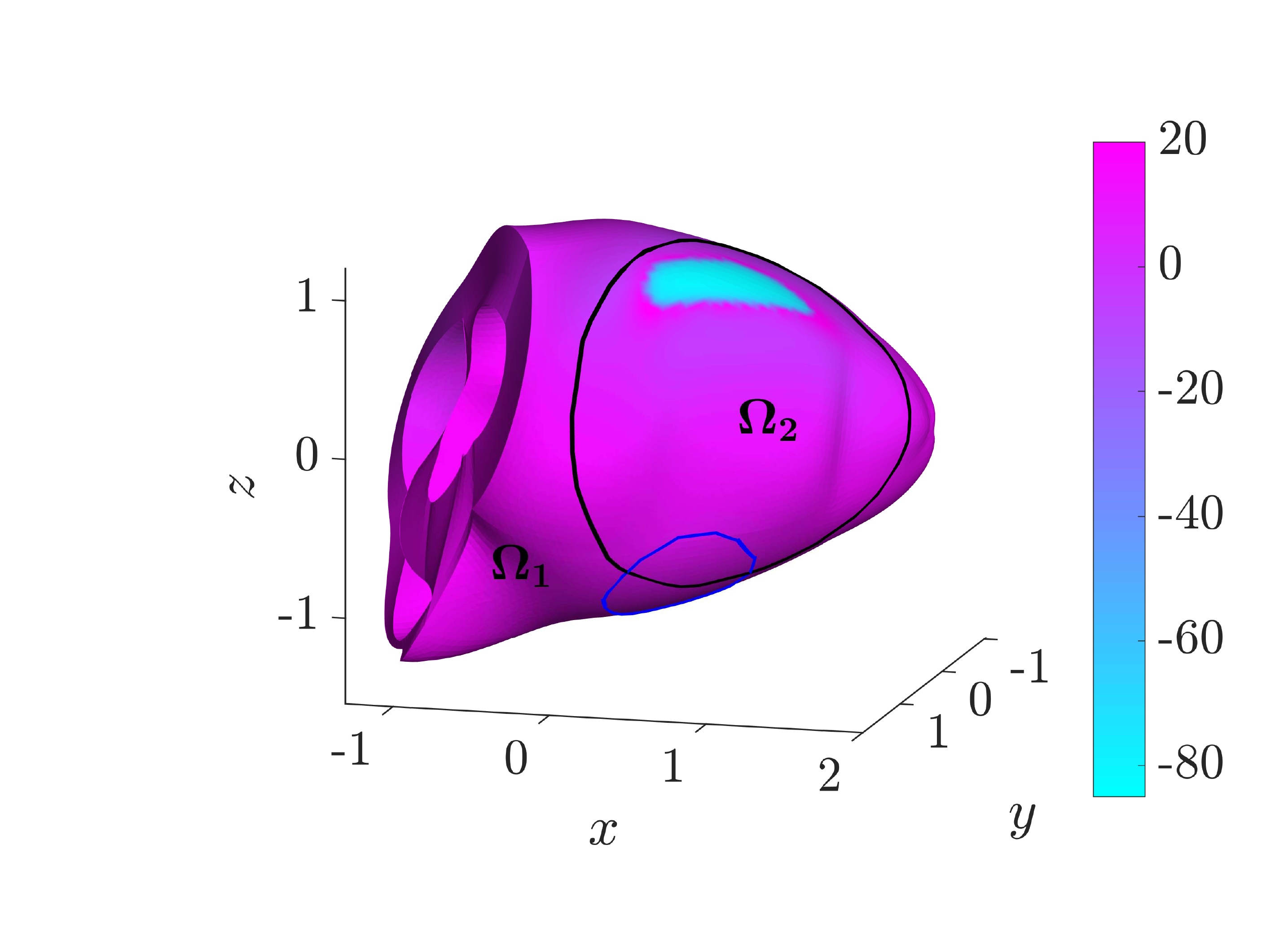}} \\
	\subfloat[$t=200$]{\includegraphics[trim=140 40 10 50, clip,width=0.49\textwidth]{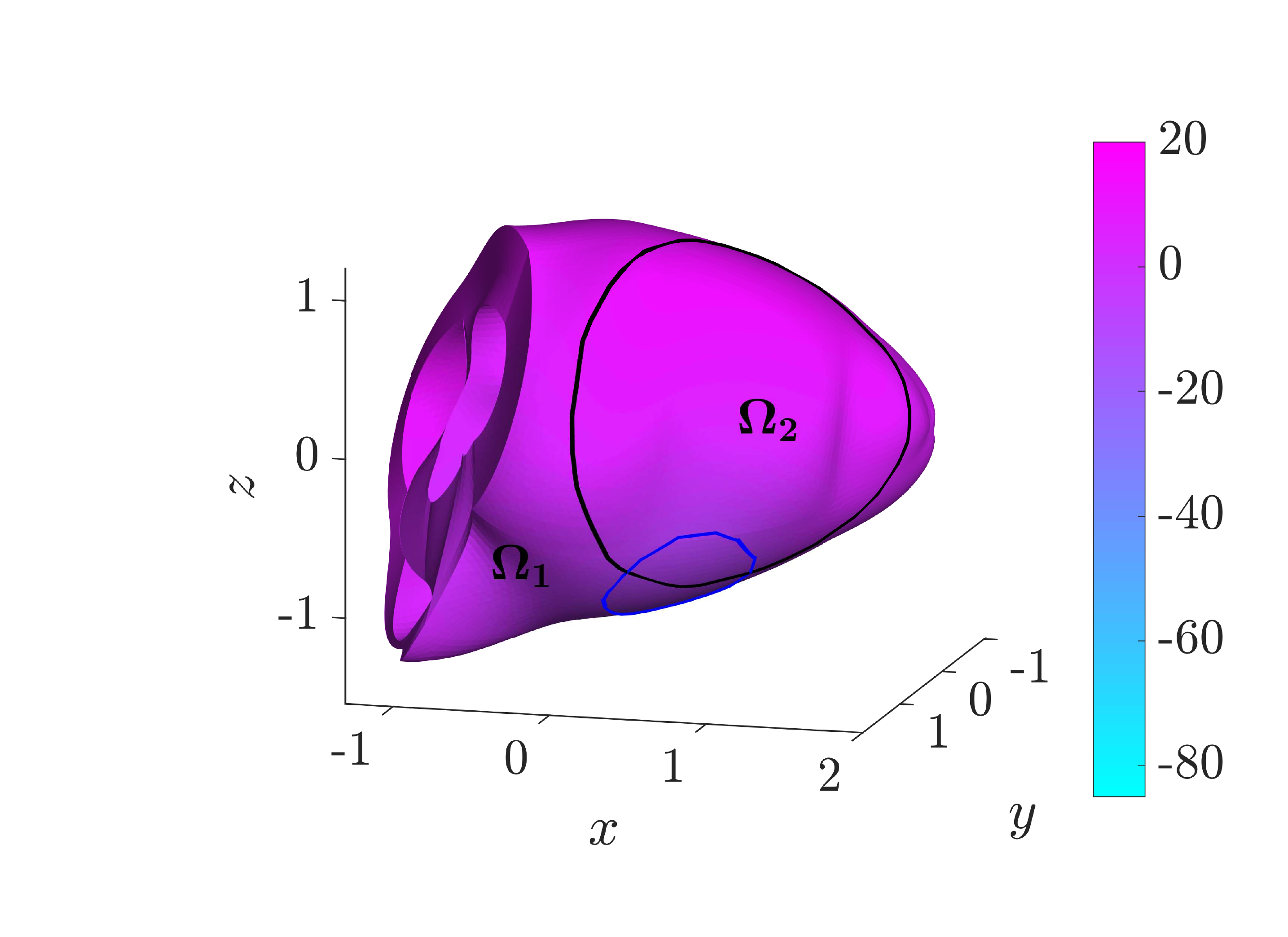}}
	\enskip
	\subfloat[$t=300$]{\includegraphics[trim=140 40 10 50, clip,width=0.49\textwidth]{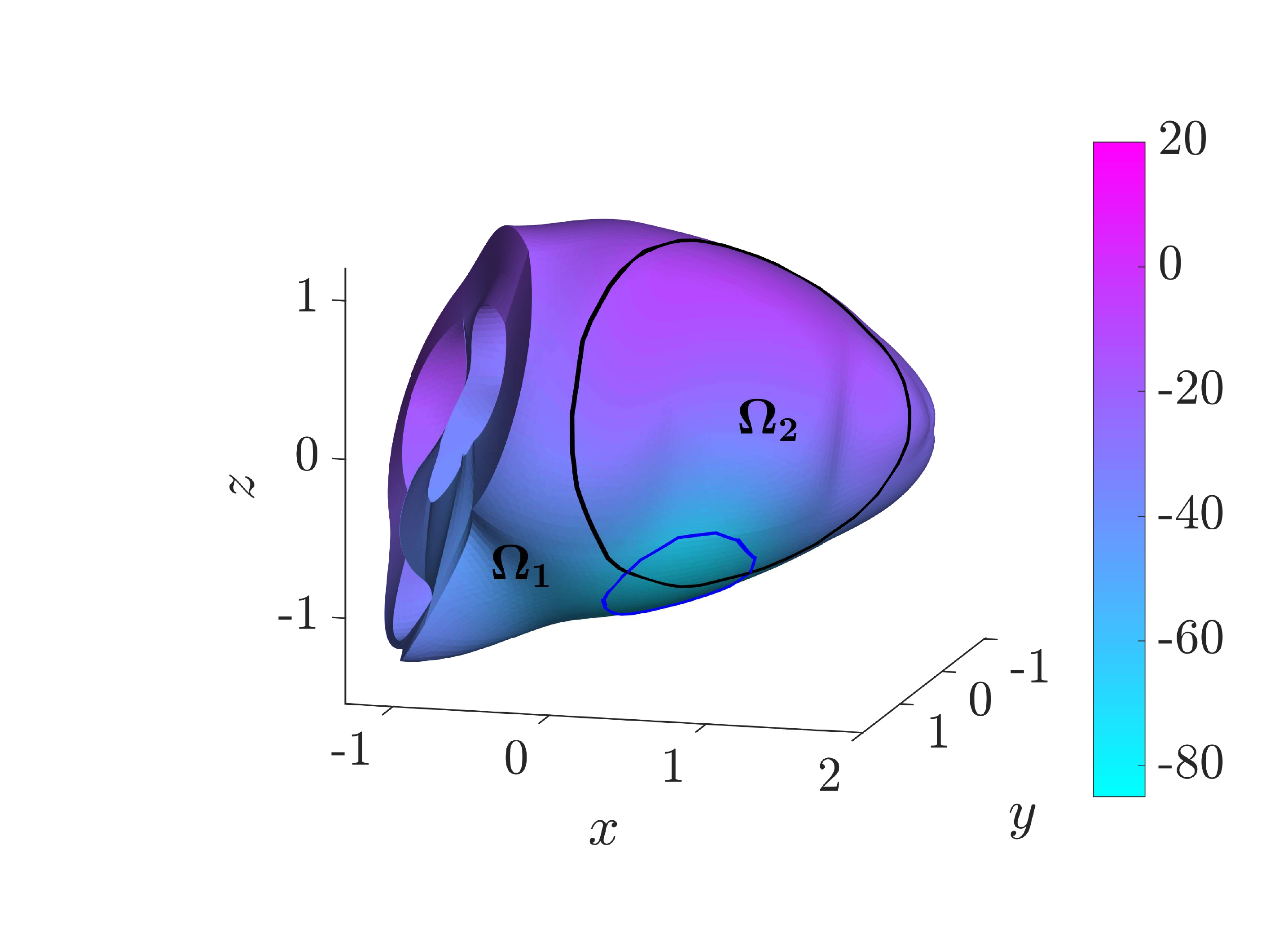}}
	\caption{Variable-order fractional monodomain problem with Beeler-Reuter Ionic current modelling re-entry with $\alpha_1=2, \alpha_2=1.7$, with stimulus applied as marked by the blue circle. Solution after the first pulse, before the second pulse.}
	\label{fig:3dresreentry1}
\end{figure}
That is, region 1 (undamaged) is modelled with a standard diffusion operator with $\alpha = 2$, while region 2 (ischaemic) uses a fractional diffusion operator with $\alpha = 1.7$.

In Figure~\ref{fig:3dresreentry1} we present the propagation of the wave after the first stimulus current is applied. In this figure we see that the tissue becomes excited at the stimulus location and the wave propagates through $\Omega_2$ at a slower rate than through $\Omega_1$.

\begin{figure}
	\centering
	\subfloat[$t=340$]{\includegraphics[trim=140 40 10 50, clip,width=0.49\textwidth]{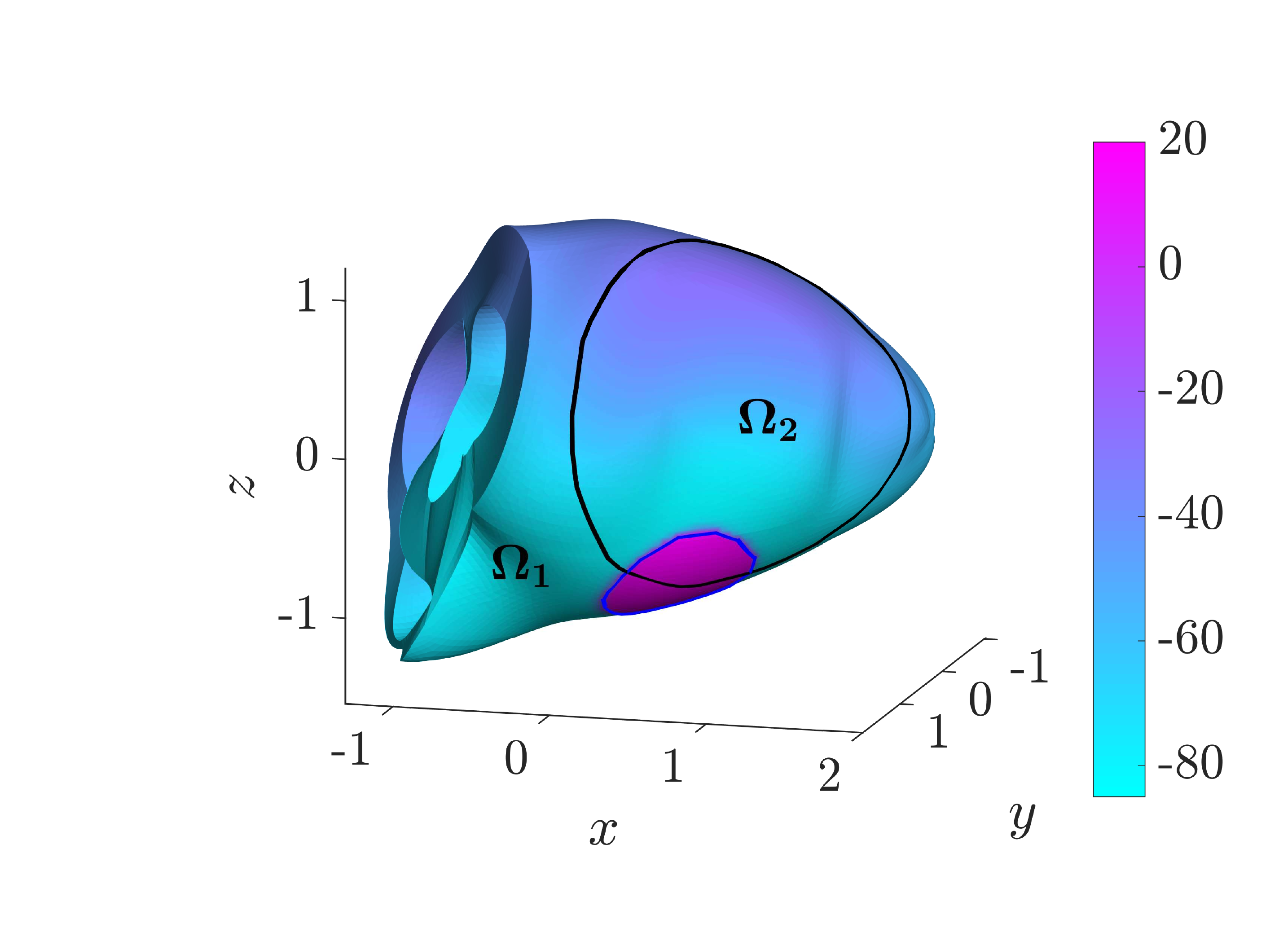}} \enskip
	\subfloat[$t=400$]{\includegraphics[trim=140 40 10 50, clip,width=0.49\textwidth]{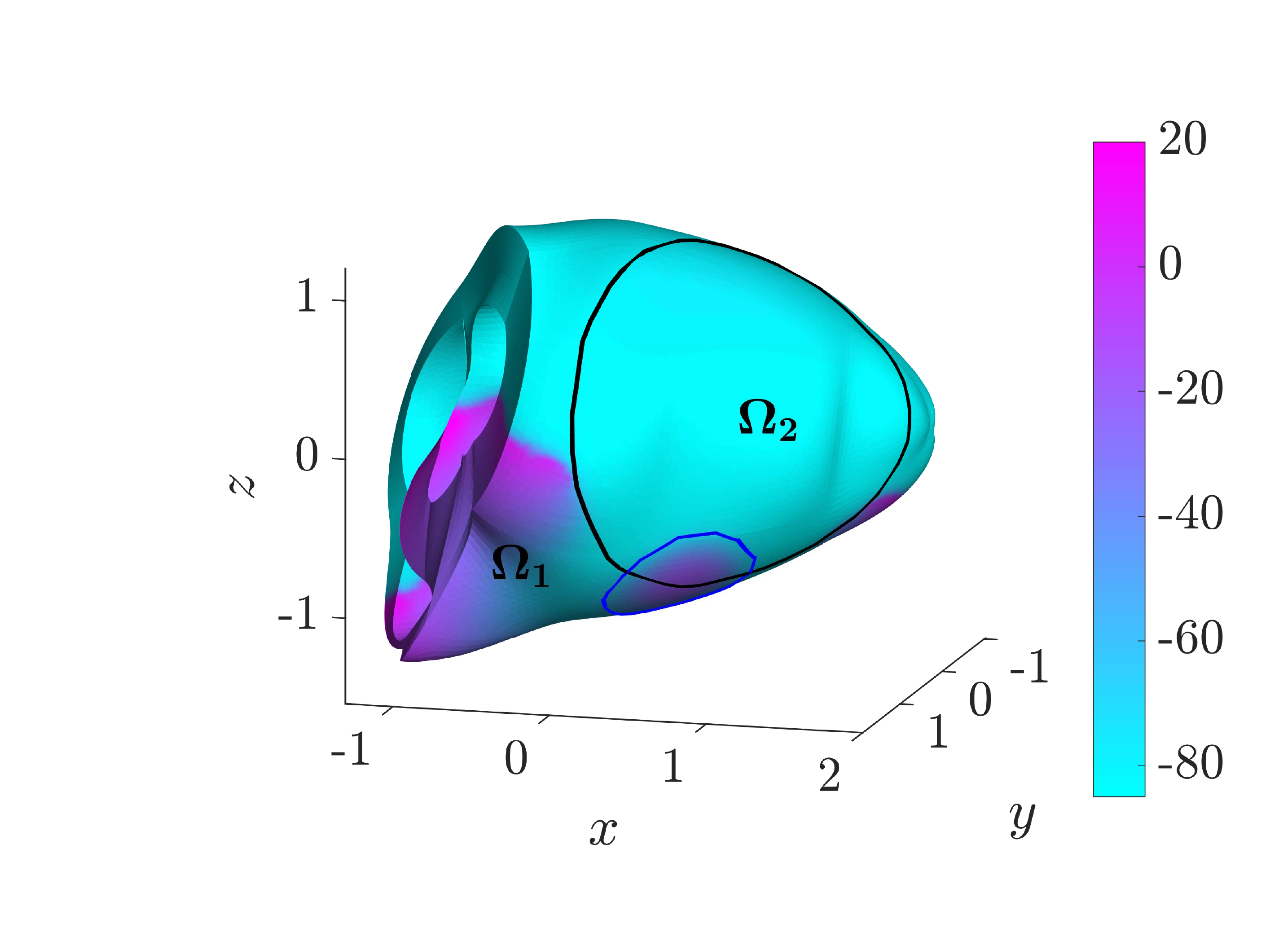}}
	\\
	\subfloat[$t=475$]{\includegraphics[trim=140 40 10 50, clip,width=0.49\textwidth]{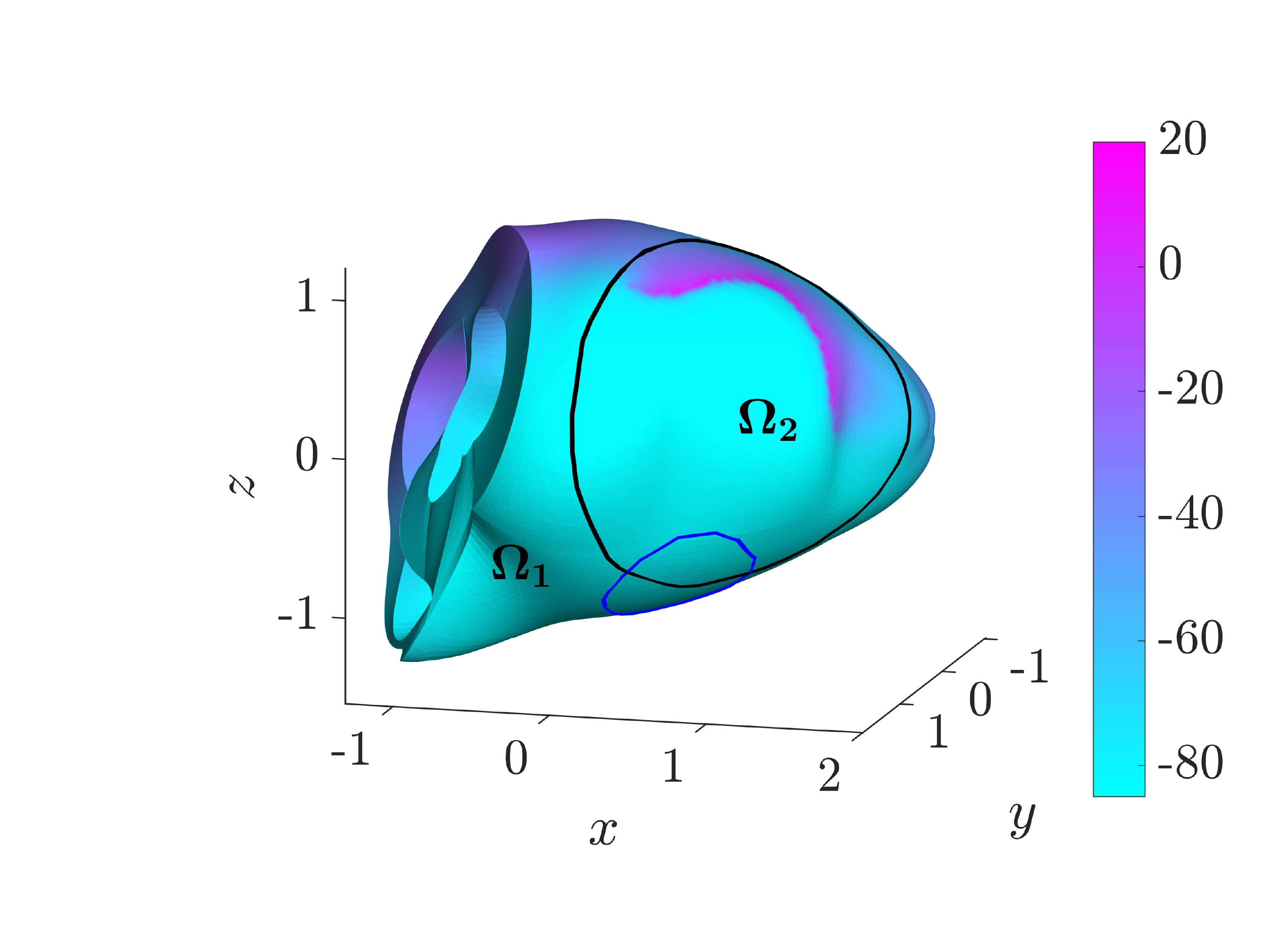}} \enskip
	\subfloat[$t=530$]{\includegraphics[trim=140 40 10 50, clip,width=0.49\textwidth]{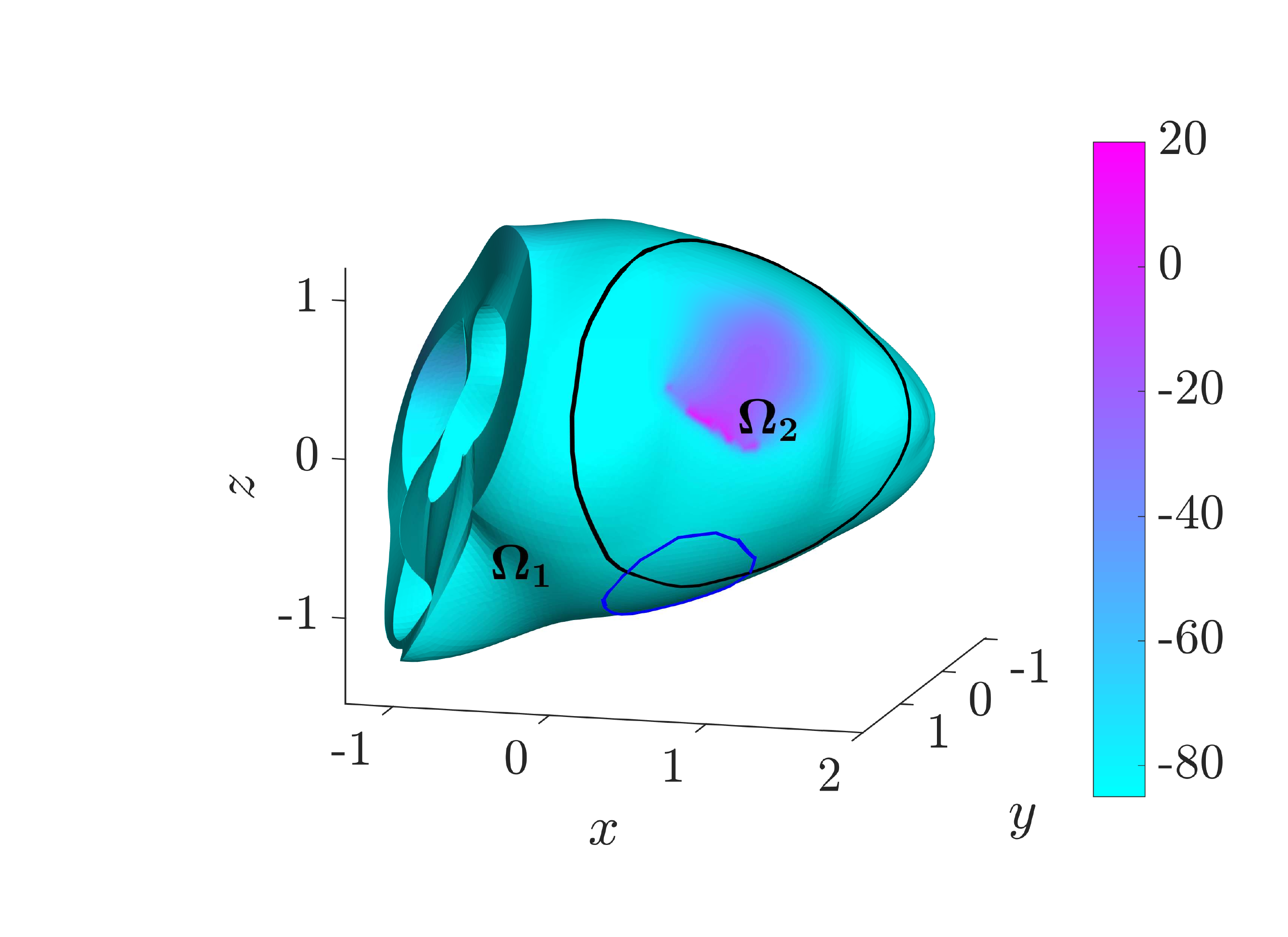}} \\
	\subfloat[$t=650$]{\includegraphics[trim=140 40 10 50, clip,width=0.49\textwidth]{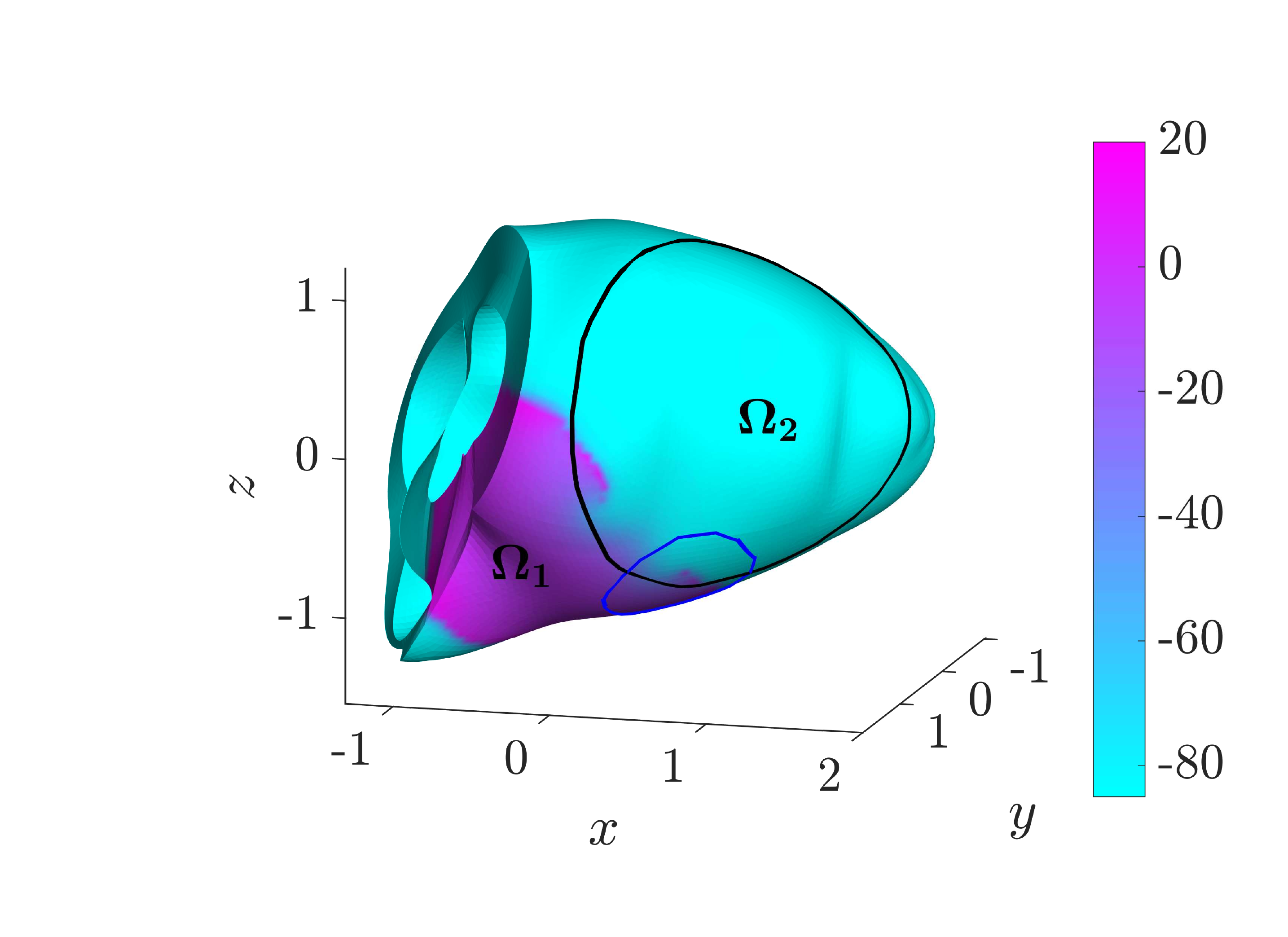}}
	\enskip
	\subfloat[$t=680$]{\includegraphics[trim=140 40 10 50, clip,width=0.49\textwidth]{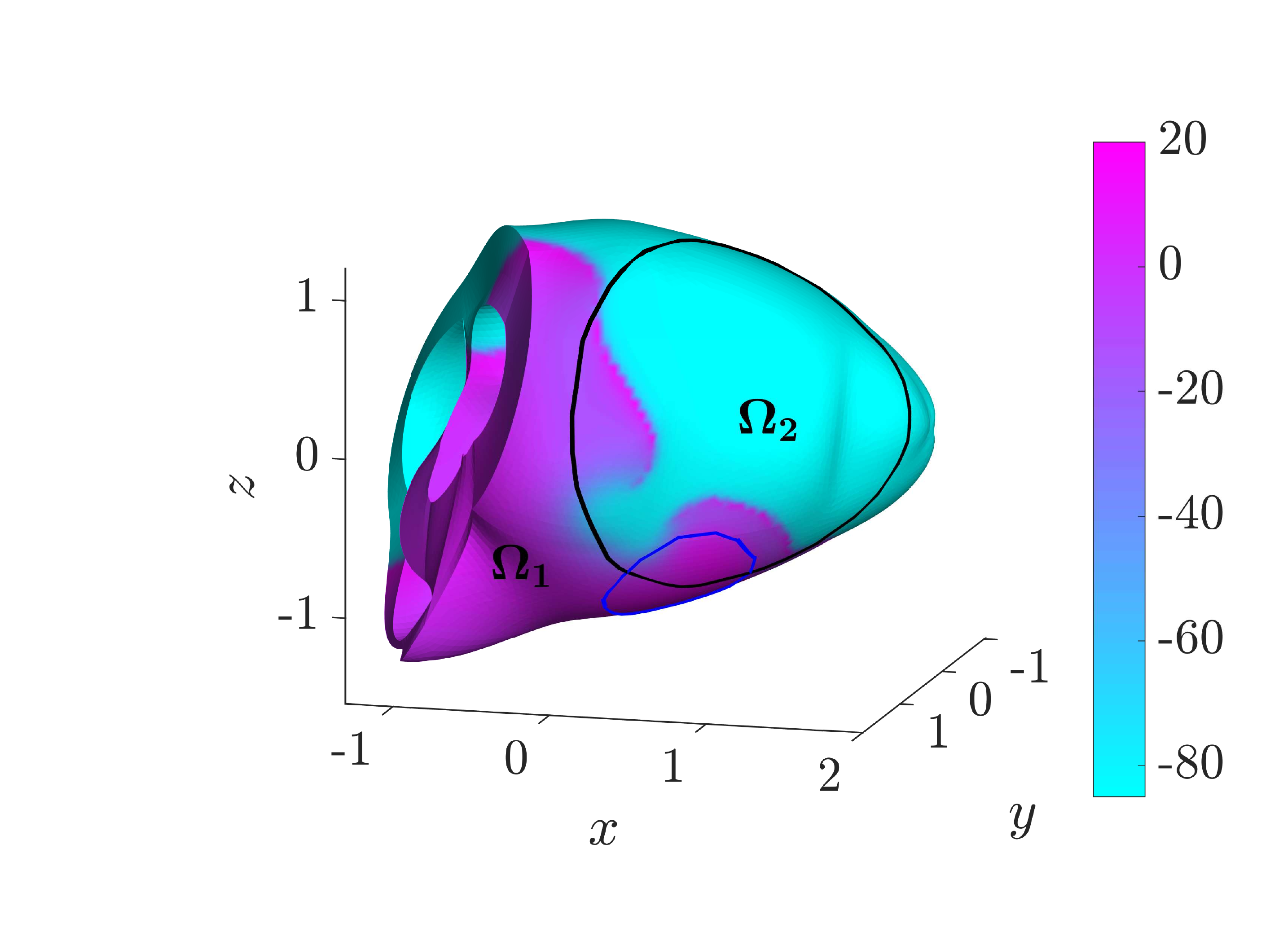}} \caption{Variable-order fractional monodomain problem with Beeler-Reuter Ionic current modelling re-entry with $\alpha_1=2, \alpha_2=1.7$, with stimulus applied as marked by the blue circle. Solution after the second pulse, before the third pulse.}
	\label{fig:3dresreentry2}
\end{figure}

In Figure~\ref{fig:3dresreentry2} we present the results showing the propagation of the electrical stimulus after the second electrical impulse has been applied. In these figures, we observe that the applied stimulus does not propagate through $\Omega_2$ in the same way. Rather, the wave travels through the domain, around $\Omega_2$, through the tissue between the ventricles and wall of the other ventricle and exhibits retrograde propagation through $\Omega_2$. Once the impulse has travelled the whole way through $\Omega_2$, it re-enters the healthy tissue, $\Omega_1$.

These results represent the culmination of the numerical and computational groundwork laid in this work.  They demonstrate the applicability of the proposed solution strategy to fully unstructured meshes on irregular domains in three dimensions, with the avoidance of dense matrices the key enabler for this advancement.

Furthermore, the results show that wave re-entry behaviour can be observed with the only modification to the standard monodomain equation being the introduction of fractional diffusion in the ischaemic region {and an appropriate choice for the time between stimuli}.  Though this work has been mostly focused on developing novel computational methods, these results suggest that accounting for spatial heterogeneities though variability in the fractional order could be an important addition in realistic mathematical models for simulating disturbances in cardiac electrical activity.

\section{Conclusions}
\label{sec:conclusions4}

In this work we developed a numerical representation of the variable-order fractional Laplacian operator where the fractional order varies across two regions representing healthy tissue and damaged tissue. The representation that we proposed can be utilised on any domain on which a matrix representation of the standard Laplacian operator can be generated. Through the use of Krylov subspace iterative methods for matrix functions, it avoids the need to compute with dense matrices.  Thus it is readily extensible to irregular meshes and higher dimensions.

We perform a number of numerical experiments to showcase the performance of our method. In comparison to previous work by Cusimano~\cite{cusimano15fractional} in one dimension, we observe complete agreement between solutions generated using our method and those of Cusimano.  Comparisons between runtime characteristics of the methods show that the iteration cost of our method approaches that of Cusimano for larger one-dimensional problems, but is absent the up-front cost of forming and factorising a dense coefficient matrix.

We further presented results in both one and three dimensions for the variable-order fractional monodomain equation with Beeler-Reuter ionic current. In these results we tested the effect of the introduction of a variable-order fractional operator.  Notably, in three dimensions we were able to observe wave re-entry behaviour, brought about only by varying the value of the fractional order in a region representing damaged tissue.  Future research will utilise this efficient computational framework to explore the implications of this observation on further understanding the factors that influence disturbances in cardiac electrical activity.

\section*{Acknowledgements}
The Authors thank Dr Nicole Cusimano for providing codes to reproduce some of the results in her PhD thesis. Author Farquhar acknowledges the support of the Australian Research Council Centre of Excellence in Mathematical and Statistical Frontiers (ACEMS).  Authors Moroney and Turner acknowledge that this research was partially supported by the Australian Research Council(ARC) via the Discovery Project (DP180103858). Author Yang acknowledges that this research was partially supported by the ARC via the Discovery Early Career Researcher Award (DE150101842).

\bibliographystyle{plain}
\bibliography{mypapers2abbrev}

\appendix
\section{Ionic current model - Beeler-Reuter}
\label{sec:BRmodel}
The Beeler-Reuter ionic current model incorporates the effects of the currents of Calcium, Sodium and Potassium. The ionic current is comprised of the sum of four different components and is written as
\begin{equation}
I_{ion} = I_{Na} + I_{K} + I_x + I_s,
\end{equation}
where $I_{Na}$ is the current carried by Sodium
\begin{equation}
I_{Na} = (4m^3 h j+0.003)(v-50),
\end{equation}
$I_{K}$ and $I_x$ are the potassium currents defined by
\begin{equation}
I_{K} = 1.4 \frac{\exp(0.04(v+85))-1}{\exp(0.08(v+53))+\exp(0.04(v+53))} + 0.07 \frac{v+23}{1-\exp(-0.04(v+23))}
\end{equation}
and
\begin{equation}
I_x = 0.8 x \frac{\exp(0.04(v+77))-1}{\exp(0.04(v+35))}
\end{equation}
and $I_s$ is the current of calcium, given by
\begin{equation}
I_s = 0.09 d \  f(v+82.3+13.0287\ln(10^{-7}c)).
\label{eq:CaCurr}
\end{equation}
The ionic current is controlled by six gating variables, $m,h,j,d,f$ and $x$, and the intracellular calcium concentration $c = 10^7\left[Ca\right]_i$, which has been scaled to simplify notation. In these equations the units of all the currents are in $\mu A\cdot cm^{-2}$, $v$ is in $mV$, the gating variables are dimensionless and $\left[Ca\right]_i$ is in moles per litre.

The equations describe the fast inward current carried by sodium, a slow inward current that is mostly carried by calcium and two outward potassium currents. The equations that describe the behaviour in a single cell consist of a system of eight ODEs
\begin{equation}
\begin{aligned}
C_m\frac{dv}{dt} &= -I_{ion}+I_{stim},\\
\frac{dG}{dt} &= \alpha_G(v)(1-G)-\beta_G(v)G, \ \ G = \{m,h,j,d,f,x\}\\
\frac{dc}{dt} &= 0.07(1-c) - I_s,
\end{aligned}
\label{eq:BRode1}
\end{equation}
where $G$ are the generic gating variables and the functions $\alpha_G$ and $\beta_G$ describe the channel opening and closing rates for the variable $G = m,h,j,d,f,x$ that the equation is referring to. For all six gating variables, both $\alpha_G$ and $\beta_G$ have a similar form
\begin{equation}
\frac{C_1\exp(C_2(v+C_3))+C_4(v+C_5)}{\exp(C_6(v+C_3))+C_7},
\end{equation}
where the values of the coefficients are different for each of the $\alpha_G$ and $\beta_G$  and can be seen in Table~\ref{tab:beelerreuter5}.

\begin{table}
	\centering
	\begin{tabular}{cccccccc}
		\toprule
		& $C_1$ & $C_2$ & $C_3$ & $C_4$ & $C_5$ & $C_6$ & $C_7$ \\
		\midrule
		$\alpha_m$ & $0$ & $0$ & $47$ & $-1$ & $47$ & $-0.1$ & $-1$\\
		$\beta_m$ & $40$ & $-0.056$ & $72$ & $0$ & $0$ & $0$ & $0$\\
		\midrule
		$\alpha_h$ & $0.126$ & $-0.25$ & $77$ & $0$ & $0$ & $0$ & $0$\\
		$\beta_h$ & $1.7$ & $0$ & $22.5$ & $0$ & $0$ & $-0.082$ & $1$\\
		\midrule
		$\alpha_j$ & $0.055$ & $-0.25$ & $78$ & $0$ & $0$ & $-0.2$ & $1$\\
		$\beta_j$ & $0.3$ & $0$ & $32$ & $0$ & $0$ & $-0.1$ & $1$\\
		\midrule
		$\alpha_d$ & $0.095$ & $-0.01$ & $-5$ & $0$ & $0$ & $-0.072$ & $1$\\
		$\beta_d$ & $0.07$ & $-0.017$ & $44$ & $0$ & $0$ & $0.05$ & $1$\\
		\midrule
		$\alpha_f$ & $0.012$ & $-0.008$ & $28$ & $0$ & $0$ & $0.15$ & $1$\\
		$\beta_f$ & $0.0065$ & $-0.02$ & $30$ & $0$ & $0$ & $-0.2$ & $1$\\
		\midrule
		$\alpha_x$ & $0.0005$ & $0.083$ & $50$ & $0$ & $0$ & $0.057$ & $1$\\
		$\beta_x$ & $0.0013$ & $-0.06$ & $20$ & $0$ & $0$ & $-0.04$ & $1$\\
		\bottomrule
		\\
	\end{tabular}
	\caption{Parameter Values for Beeler-Reuter model\cite{beeler77reconstruction}}
	\label{tab:beelerreuter5}
\end{table}

When the Beeler-Reuter model is coupled with the monodomain equation to include spatial dependence of the potential, the result is the coupled PDE-ODE system
\begin{equation}
\begin{aligned}
\chi \left(C_m\frac{d v}{dt} + I_{ion}\right) &= -\frac{\lambda}{1+\lambda}M_i \left(-\nabla^2  \right) v+ I_{stim},\\
\frac{dG}{dt} &= \alpha_G(v)(1-G)-\beta_G(v)G,\\
\frac{dc}{dt} &= 0.07(1-c) - I_s.
\end{aligned}
\label{eq:BRmodel}
\end{equation}
These equations can be used to describe how electrical impulses propagate through a spatial domain, such as the heart.
\end{document}